\documentclass[11pt]{article} \topmargin=-0.5cm
\usepackage{amssymb}

   \csname @addtoreset\endcsname{equation}{section}
    \csname @addtoreset\endcsname{figure}{section}
    \csname @addtoreset\endcsname{table}{section}

\newcounter{thanksnum}
\def\thanksnumber#1
{\setcounter{thanksnum}{\value{footnote}}\setcounter{footnote}{#1}%
                     \addtocounter{footnote}{-1}\footnotemark
                     \setcounter{footnote}{\value{thanksnum}}}
\def\newtheoremz#1{\@ifnextchar[{\@othmz{#1}}{\@nthmz{#1}}}

\def\@nthmz#1#2{%
\@ifnextchar[{\@xnthmz{#1}{#2}}{\@ynthmz{#1}{#2}}}

\def\@xnthmz#1#2[#3]{\expandafter\@ifdefinable\csname #1\endcsname
{\@definecounter{#1}\@addtoreset{#1}{#3}%
\expandafter\xdef\csname the#1\endcsname{\expandafter\noexpand
  \csname the#3\endcsname \@thmcountersepz \@thmcounterz{#1}}%
\global\@namedef{#1}{\@thmz{#1}{#2}}\global\@namedef{end#1}{\@endtheoremz}}}

\def\@ynthmz#1#2{\expandafter\@ifdefinable\csname #1\endcsname
{\@definecounter{#1}%
\expandafter\xdef\csname the#1\endcsname{\@thmcounterz{#1}}%
\global\@namedef{#1}{\@thm{#1}{#2}}\global\@namedef{end#1}{\@endtheoremz}}}

\def\@othmz#1[#2]#3{\expandafter\@ifdefinable\csname #1\endcsname
  {\global\@namedef{the#1}{\@nameuse{the#2}}%
\global\@namedef{#1}{\@thmz{#2}{#3}}%
\global\@namedef{end#1}{\@endtheoremz}}}

\def\@thmz#1#2{\refstepcounter
    {#1}\@ifnextchar[{\@ythmz{#1}{#2}}{\@xthmz{#1}{#2}}}

\def\@xthmz#1#2{\@begintheoremz{#2}{\csname the#1\endcsname}\ignorespaces}
\def\@ythmz#1#2[#3]{\@opargbegintheoremz{#2}{\csname
       the#1\endcsname}{#3}\ignorespaces}

\def\@thmcounterz#1{\noexpand\arabic{#1}}
\def\@thmcountersepz{.}
\def\@begintheoremz#1#2{ \trivlist \item[\hskip \labelsep{\bf #1\ #2}]}
\def\@opargbegintheoremz#1#2#3{ \trivlist
      \item[\hskip \labelsep{\bf #1\ #2\ (#3)}]}
\def\@endtheoremz{\endtrivlist}

\newtheorem{theorem}{Theorem}[section]
\newtheorem{lemma}{Lemma}[section]

\newtheorem{proposition}{Proposition}[section]
\newtheorem{corollary}{Corollary}[section]
\newtheorem{condition}{Condition}[section]
\newtheorem{definition}{Definition}[section]
\newtheorem{remark}{Remark}[section]

\newtheorem{example}{Example}[section]

\def\e{\varepsilon}

\def\defi{\stackrel{{\scriptscriptstyle \Delta}}{=}}

\def\a{\alpha}
\def\d{\delta}
\def\o{\omega}
\def\O{\Omega}

\def\Y{{\cal Y}}
\def\F{{\cal F}}
\def\w{\widehat}
\def\Ind{{\mathbb{I}}}

\def\esssup{\mathop{\rm ess\, sup}}
\def\const{{\rm const\,}}

\def\R{{\bf R}}
\def\E{{\bf E}}
\def\P{{\bf P}}

\def\L{L}

\def\b{\beta}
\def\s{\delta}
\def\g{\gamma}

\def\W{{\cal W}^*}
\def\ww{\widetilde}
\def\X{{\cal X}}
\def\t{\theta}
\def\oo{\bar}
\def\s{\sigma}

\def\p{\partial}
\def\G{\Gamma}

\def\A{{\cal A}}
\def\M{{\cal M}}

\def\L{{\cal L}}
\def\I{{\, \cal I}}

\newcommand{\be}{\begin{equation}}
\newcommand{\ee}{\end{equation}}
\newcommand{\bd}{\begin{displaymath}}
\newcommand{\ed}{\end{displaymath}}
\newcommand{\ba}{\begin{array}{ll}}
\newcommand{\ea}{\end{array}}
\newcommand{\baa}{\begin{eqnarray}}
\newcommand{\eaa}{\end{eqnarray}}
\newcommand{\baaa}{\begin{eqnarray*}}
\newcommand{\eaaa}{\end{eqnarray*}}   \font\sm=cmr10

\def\ww{\tilde}
\def\W{{\cal W}}

\def\Q{{\cal Q}}
\def\CC{{\cal C}}

\def\LU{L}

\def\KK{{K}}

\def\I{{ \d}}

\date{First version: 23 June 2006. Revision:  26 July 2010}
\title{Backward parabolic Ito equations and second fundamental inequality}
\author{ Nikolai Dokuchaev
\\ {\sm  Department of Mathematics and Statistics,
Curtin University of Technology},\\ {\sm GPO Box U1987, Perth,
Western Australia, 6845}}
\begin{document}
\maketitle
\begin{abstract}
Regularity of solutions is studied for backward stochastic
parabolic Ito equations. An analog of the second energy inequality
and the related existence theorem are obtained for domains with
boundary.
 \\ {\it Key
words and phrases:} parabolic Ito equations, backward SPDEs,
regularity
\end{abstract}
\section{Introduction}
The paper studies stochastic partial differential equations (SPDEs)
in a cylinder $D\times[0,T]$ with a Dirichlet boundary condition on
$\p D$, for a region $D\subseteq \R^n$. We investigate regularity
properties of the backward equations, i.e., equations  with Cauchy
condition at the final time. The difference between backward and
forward equations is not that important for the deterministic
equations since a deterministic backward equation can be converted
to a forward equation by a time change. It cannot be done so easily
for stochastic equations, because we look for solutions adapted to
the driving Brownian motion.
 It is why the backward
SPDEs require special consideration. The most common  approach is
to consider the so-called Bismut backward equations such the
diffusion term is not given a priori but needs to be found. These
approach was introduced first for ordinary linear backward
stochastic equations. The backward SPDEs with similar features
were widely studied (see, e.g., Pardoux and Peng
 (1990), Hu and Peng  (1991),  Dokuchaev (1992),
 (2003),
  Yong and Zhou
 (1999),
 Pardoux and Rascanu (1998), Ma and Yong (1999), Hu {\it et al} (2002), Confortola
 (2007), and the bibliography there given). Backward parabolic SPDEs represent analogs of backward parabolic
Kolmogorov equations for non-Markov Ito processes, including the
case of bounded domains, so they can be used for characterization of
distributions of the first exit times in non-Markovian setting, as
was shown by the author (1992,2010a).  A different type of backward
equations was described in Chapter 5 of Rozovskii (1990).
  Forward SPDEs were also widely studied (see,
e.g., Al\'os et al (1999), Bally {\it et al} (1994),
Chojnowska-Michalik and Goldys (1995), Da Prato and Tubaro (1996),
Gy\"ongy (1998), Krylov (1999), Maslowski (1995), Pardoux (1993),
 Rozovskii (1990), Walsh (1986), Zhou (1992), Dokuchaev
 (1995), (2002), (2005),
and the bibliography there given).
\par
For linear PDEs,  existence and uniqueness at different spaces is
expressed traditionally via a priori estimates, when a norm of the
solution is estimated via a norm of the free term. For the second
order equations, there are two most important estimates based on the
$L_2$-norm: so-called "the first energy inequality" or "the first
fundamental inequality", and "the second energy inequality", or "the
second fundamental inequality" (Ladyzhenskaya (1985)). For instance,
consider a boundary value problem for the  heat equation \baa
&&u_t'=u''_{xx}+\varphi,\quad\varphi=f'_x+g,\nonumber
\\&&u|_{t=0}=0,\quad u|_{\p D}=0,\quad (x,t)\in Q=D\times [0,1], \quad
D\subset \R. \label{simple} \eaa Then the first fundamental
inequality is the estimate $$
\|u'_x\|^2_{L_2(Q)}+\|u\|^2_{L_2(Q)}\le\const
(\|f\|^2_{L_2(Q)}+\|g\|^2_{L_2(Q)}).$$ Respectively, the second
fundamental inequality is the estimate $$
\|u\|^2_{L_2(Q)}+\|u'_x\|^2_{L_2(Q)}+\|u''_{xx}\|^2_{L_2(Q)}
\le\const \|\varphi\|^2_{L_2(Q)}.$$
  Note that the
second fundamental inequality leads to existence theorem in the
class of solutions such that $u''_{xx}\in L_2(Q)$, and the first
fundamental inequality leads to existence theorem in the class of
solutions such that $u'_x\in L_2(Q)$, i.e., with generalized
distributional derivatives $u''_{xx}$ only.  For the problem without
boundary value condition, with smooth coefficients, and
one-dimensional  $x\in \R$, the second fundamental inequality can be
derived from the first fundamental inequality; it suffices to apply
the first fundamental inequality for the parabolic equation for
$u'_x$. (For the vector case of $x\in\R^n$, it would be more
difficult since $u'_x$ is a vector satisfying a system of $n$
parabolic equations). For the problems with boundary value
conditions, this approach does not work even for one-dimensional
case, since the boundary values on $\p D$ for $u'_x$ are unknown a
priori. It is why the second fundamental inequality needs to be
derived separately using special methods.
\par
For forward parabolic SPDEs, analogs of the first and the second
fundamental inequalities are known. These results are summarized in
Lemma \ref{lemma1}   below. The first fundamental inequality for
forward SPDEs in bounded domains with Dirichlet boundary condition
was known long time ago (see, e.g.,  Rozovskii (1990)). Moreover,
similar results are also known for forward SPDEs of an arbitrary
high order $2m\ge 2$; in this setting, the analog of "the first
fundamental inequality" is an estimate for
$\E\|u(\cdot,t)\|^2_{W_2^m(D)}$ (Rozovskii (1991)). In addition, a
priori estimates  without Dirichlet conditions, i.e., in the  entire
space, are known  for a general setting that covers both first and
second fundamental inequalities (Krylov (1999)). On the other hand,
"the second fundamental inequality" for the problem with boundary
conditions was more difficult to obtain. Related complications were
discussed in Krylov (1999), p. 237 and in Dokuchaev (2005). Kim
(2004) obtained a priori estimates for forward parabolic SPDEs for
special weighted norms  that devaluates boundary values; for the
case of $L_2$-norms, these estimates can be interpreted as analogs
of "the second fundamental inequality"; they are similar to
estimates $\|r_1u''_{xx}\|_{L_2(Q)}\le\const
\|r_2\varphi\|_{L_2(Q)}$ for the problem (\ref{simple}), where $r_i$
are some weight functions such that $r_i(x)\to 0$ as $x$ approaching
$\p D$. For the standard non-weighted Sobolev norms, the second
fundamental inequality"  was obtained in Dokuchaev (2005).
\par
For the backward parabolic equations with Dirichlet boundary
conditions, an analog of the first fundamental inequality is known
(Zhou (1992), Dokuchaev (1992), (2003)). In fact, the duality
relationship
 between forward and backward equations makes
  it sufficient to
prove the first fundamental inequality for any one type of these two
types of equations. (By duality we mean equations (\ref{duality})
connecting the solutions of SPDEs (\ref{4.1}) and (\ref{4.8})
respectively). However, this approach does not work for the second
fundamental inequality in a bounded domain $D$, since it requires to
study an adjoint equation with the free term taking values in the
space $(W_2^2(D))^*$ which is too wide. It was unknown if the second
fundamental inequality holds in this case.
\par In the present paper, we study again
existence, uniqueness, and a priori estimates for solutions for
backward SPDEs. As was mentioned above, the first and the second
fundamental inequalities for the forward SPDEs had been proved, as
well as the first fundamental inequality for the backward SPDEs, so
we concentrate our efforts on the remaining problem: to investigate
if an analog of the second fundamental inequality holds for the
backward equations. We found  sufficient conditions that ensure that
the second fundamental inequality and the related existence theorem
holds (Theorem \ref{Th2FE}).  To ensure this regularity, we required
additional Condition \ref{condA} which is  a strengthened version of
the standard coercivity condition (Condition \ref{cond3.1.A}).
Without this new condition, the second fundamental inequality is
still not established.
\par
Some examples of applications are discussed in Section
\ref{Secappl}.
\section{Definitions}\label{SecD}
\subsection{Spaces and classes of functions.} 
Assume that we are given an open domain $D\subseteq\R^n$ such that
either $D=\R^n$ or $D$ is bounded  with $C^2$-smooth boundary $\p
D$. Let $T>0$ be given, and let $Q\defi D\times (0,T)$. \par We
are given a standard complete probability space $(\O,\F,\P)$ and a
right-continuous filtration $\F_t$ of complete $\s$-algebras of
events, $t\ge 0$. We are given also a $N$-dimensional process
$w(t)=(w_1(t),...,w_N(t))$ with independent components such that
it is a Wiener process with respect to $\F_t$.
\par
We denote by $\|\cdot\|_{ X}$ the norm in a linear normed space
$X$, and
 $(\cdot, \cdot )_{ X}$ denote  the scalar product in  a Hilbert space $
X$.
\par
We introduce some spaces of real valued functions.
\par
 Let $G\subset \R^k$ be an open
domain, then ${W_q^m}(G)$ denote  the Sobolev  space of functions
that belong to $L_q(G)$ with the distributional derivatives up to
 the $m$th order, $q\ge 1$.
\par
 We denote by $|\cdot|$ the Euclidean norm in $\R^k$, and we denote by $\bar G$
the closure of a region $G\subset\R^k$.
\par Let $H^0\defi L_2(D)$,
and let $H^1\defi \stackrel{\scriptscriptstyle 0}{W_2^1}(D)$ be
the closure in the ${W}_2^1(D)$-norm of the set of all smooth
functions $u:D\to\R$ such that  $u|_{\p D}\equiv 0$. Let
$H^2=W^2_2(D)\cap H^1$ be the space equipped with the norm of
$W_2^2(D)$. The spaces $H^k$ are Hilbert spaces, and $H^k$ is a
closed subspace of $W_2^k(D)$,  $k=1,2$.
\par
 Let $H^{-1}$ be the dual space to $H^{1}$, with the
norm $\| \,\cdot\,\| _{H^{-1}}$ such that if $u \in H^{0}$ then
$\| u\|_{ H^{-1}}$ is the supremum of $(u,v)_{H^0}$ over all $v
\in H^1$ such that $\| v\|_{H^1} \le 1 $. $H^{-1}$ is a Hilbert
space.
\par We will write $(u,v)_{H^0}$ for $u\in H^{-1}$
and $v\in H^1$, meaning the obvious extension of the bilinear form
from $u\in H^{0}$ and $v\in H^1$.
\par
We denote by $\oo\ell _{k}$ the Lebesgue measure in $\R^k$, and we
denote by $ \oo{{\cal B}}_{k}$ the $\sigma$-algebra of Lebesgue
sets in $\R^k$.
\par
We denote by $\oo{{\cal P}}$  the completion (with respect to the
measure $\oo\ell_1\times\P$) of the $\s$-algebra of subsets of
$[0,T]\times\O$, generated by functions that are progressively
measurable with respect to $\F_t$.
\par
Let $Q_s\defi D\times [s,T]$. For $k=-1,0,1,2$, we  introduce the
spaces
 \baaa
 &&X^{k}(s,T)\defi L^{2}\bigl([ s,T ]\times\Omega,
{\oo{\cal P }},\oo\ell_{1}\times\P;  H^{k}\bigr),\, \quad Z^k_t
\defi L^2\bigl(\Omega,{\cal F}_t,\P;
H^k\bigr),\quad \\&&\CC^{k}(s,T)\defi C\left([s,T]; Z^k_T\right).
\eaaa
Furthermore, we introduce the spaces $$ Y^{k}(s,T)\defi
X^{k}(s,T)\!\cap \CC^{k-1}(s,T), \quad k\ge 0, $$ with the norm $
\| u\| _{Y^k(s,T)}
\defi \| u\|_{{X}^k(s,T)} +\| u\| _{\CC^{k-1}(s,T)}. $
\par
In addition, we will be using spaces \baaa &&\W^{k}_{r} \defi
L^{\infty}\bigl([0,T ]\times\O, \overline{{\cal
P}},\oo\ell_{1}\times\P;\, W_r^{k}(D)\bigr), \quad k=0,1,\ldots,
\quad 1\le r\le
 +\infty.
\eaaa
\par
  The
spaces $X^k$ and $Z_t^k$  are Hilbert spaces.
\begin{proposition} 
\label{propL} Let $\xi\in X^0$,
 let a sequence  $\{\xi_k\}_{k=1}^{+\infty}\subset
L^{\infty}([0,T]\times\O, \ell_1\times\P;\,C(D))$ be such that all
$\xi_k(\cdot,t,\o)$ are progressively measurable with respect to
$\F_t$, and let $\|\xi-\xi_k\|_{X^0}\to 0$. Let $t\in [0,T]$ and
$j\in\{1,\ldots, N\}$ be given.
 Then the sequence of the
integrals $\int_0^t\xi_k(x,s,\o)\,dw_j(s)$ converges in $Z_t^0$ as
$k\to\infty$, and its limit depends on $\xi$, but does not depend
on $\{\xi_k\}$.
\end{proposition}
\par
{\it Proof} follows from completeness of  $X^0$ and from the
equality
\begin{eqnarray*}
\E\int_0^t\|\xi_{k}(\cdot,s,\o)-\xi_m(\cdot,s,\o)\|_{H^0}^2\,ds
=\int_D\,dx\,\E\left(\int_0^t\big(\xi_k(x,s,\o)-
\xi_m(x,s,\o)\big)\,dw_j(s)\right)^2.
\end{eqnarray*}
\begin{definition} 
\rm Let $\xi\in X^0$, $t\in [0,T]$, $j\in\{1,\ldots, N\}$, then we
define $\int_0^t\xi(x,s,\o)\,dw_j(s)$ as the limit  in $Z_t^0$ as
$k\to\infty$ of a sequence $\int_0^t\xi_k(x,s,\o)\,dw_j(s)$, where
the sequence $\{\xi_k\}$ is such  as in Proposition \ref{propL}.
\end{definition}
 Sometimes we will omit
$\o$.
\section{Review of existence theorems for forward and backward SPDEs}
\label{SecC}  Let
 $(x,t)\in Q$,   $\o\in\O$.
 \par
Consider the functions $b(x,t,\o):
\R^n\times[0,T]\times\O\to\R^{n\times n}$, $f(x,t,\o):
\R^n\times[0,T]\times\O\to\R^n$, $\lambda(x,t,\o):
\R^n\times[0,T]\times\O\to\R$  $\b_j(x,t,\o):
\R^n\times[0,T]\times\O\to\R^n$, $\oo\b_i(x,t,\o):$
$\R^n\times[0,T]\times\O\to\R$ that are  progressively measurable
for any $x\in \R^n$ with respect to $\F_t$.
\par
  Consider  differential operators defined
 on functions $v:D\to\R$
  \baaa\label{A}&&\A v =\sum_{i,j=1}^n\frac{\p ^2
}{\p x_i \p x_j}\, \biggl(b_{ij}(x,t,\o)\,v(x)\biggr) -\,
\sum_{i=1}^n \frac{\p }{\p x_i }\,\big( f_{i}(x,t,\o)\,v(x)\big)
+\lambda(x,t,\o)\,v(x),
\nonumber\\
 &&B_kv=-\sum_{i=1}^n \frac{\p }{\p x_i }\,\big(\beta_{k}(x,t,\o)\,v(x))+\oo \beta_k(x,t,\o)\,v(x),\qquad k=1,\ldots ,N. \eaaa Here $b_{ij}, f_i,
x_i$ are the components of $b$,$f$, and $x$.
\par
Further, consider the operators being formally adjoint to the
operators $\A$ and $B_i$: \baa &&\A^* v\defi \sum_{i,j=1}^n
b_{ij}(x,t,\o)\frac{\p^2v}{\p x_i \p x_j}(x) +\sum_{i=1}^n
f_{i}(x,t,\o)\frac{\p v}{\p x_i }(x) +\,\lambda(x,t,\o)v(x), \nonumber\\
\label{B}&& B_k^*v\defi\frac{dv}{dx}\,(x)\,\beta_k(x,t,\o) +\oo
\beta_k(x,t,\o)\,v(x),\qquad k=1,\ldots ,N. \label{AB*}\eaa
 To proceed further, we assume that Conditions
\ref{cond3.1.A}-\ref{cond3.1.B} remain in force throughout this
paper.
 \begin{condition} \label{cond3.1.A}  (Coercivity) The matrix  $b=b^\top$ is
symmetric,  bounded, and progressively measurable with respect to
$\F_t$ for all $x$, and there exists a constant $\d_1>0$ such that
\baaa
 \label{Main1} y^\top  b
(x,t,\o)\,y-\frac{1}{2}\sum_{i=1}^N |y^\top\b_i(x,t,\o)|^2 \ge \d
|y|^2 \quad\forall\, y\in \R^n,\ (x,t)\in  D\times [0,T],\ \o\in\O.
\eaaa
\end{condition}
\begin{condition}\label{cond3.1.B}
The functions
 $\lambda (x,t,\o):\R^n$  and $\oo\b_i(x,t,\o)$  are
bounded. The functions $b(x,t,\o):\R^n \times \R\times\O\to \R^{n
\times n}$, $f(x,t,\o):\R^n \times \R\times\O\to \R^n$, $\lambda
(x,t,\o):\R^n \times \R\times\O\to \R$, $\b_i(x,t,\o)$ and
$\oo\b_i(x,t,\o)$  are are differentiable in $x$ and bounded in
$(x,t,\o)$,   and \baaa \esssup_{x,t,\o} \biggl( \Bigl| \frac{\p
b}{\p x}(x,t,\o)\Bigr| + \Bigl| \frac{\p f}{\p x}(x,t,\o)\Bigr| +
|\frac{\p \b_i}{\p x}(x,t,\o)|<+\infty, \quad i=1,\ldots ,N.  \eaaa
\end{condition}
\par
We introduce the set of parameters   $$ \ba {\cal P}_1
\defi \biggl( n,\,\, D,\,\,  T\,\, \delta,\,\,\,\,
\esssup_{x,t,\o}\Bigl[| b(x,t,\o)|+ |f(x,t,\o)|+ \Bigl|\frac{\p
b}{\p x}(x,t,\o)\Bigr|+ \Bigl|\frac{\p f}{\p
x}(x,t,\o)\Bigr|\Bigr],\\ \hspace{5cm}
 \esssup_{x,t,\o,i}\Bigl[|
\b_i(x,t,\o)|+ |\oo\b_i(x,t,\o)|+ \Bigl|\frac{\p \b_i}{\p
x}(x,t,\o)\Bigr|\Bigr] \biggr). \ea
$$
\subsubsection*{Boundary value problems for forward and backward equations}
Let $s\in [0,T)$, $\varphi\in X^{-1}$, $h_i\in X^0$, and
$\Phi,\Psi\in Z^0_s$. Consider the  boundary value problem in
$D\times[s,T]$ \baa &&d_tu=\left( \A u+ \varphi\right)dt
+ \sum_{i=1}^N(B_iu+h_i)dw_i(t),\quad \quad t >s,\nonumber\\
&&u|_{t=s}=\Phi,\quad u(x,t,\o)|_{x\in \p D}=0 .\label{4.1}
 \eaa
The corresponding SPDE is a forward equation. Here
 $u=u(x,t,\o)$,
 $(x,t)\in Q$,   $\o\in\O$.
 \par
Inequality (\ref{Main1})  means that equation (\ref{4.1}) is
coercive or {\it superparabolic}, in the terminology of Rozovskii
(1990).
\par
Further, let $\xi\in X^{-1}$, and $\Psi\in Z^0_s$. Consider the
boundary value problem in $Q$ \baa  &&d_tp+
\Bigl(\A^*p+\sum_{i=1}^NB_i^*\chi_i+\xi\Bigr)\,dt=
\sum_{i=1}^N\chi_i \,dw_i(t),\quad t<T,\nonumber
\\
&&p|_{t=T}=\Psi, \quad p(x,t,\o)\,|_{x\in \p D}=0.
\label{4.8}
\eaa The corresponding SPDE is a backward equation. Here
 $p=p(x,t,\o)$, $\chi_i=\chi_i(x,t,\o)$,
 $(x,t)\in Q$,   $\o\in\O$.

\subsubsection*{The definition of solution} 
\begin{definition} 
\label{defsolltion} \rm Let $h_i\in X^0$ and $\varphi\in X^{-1}$. We
say that   equations (\ref{4.1}) are satisfied for $u\in Y^1$ if
\baa u(\cdot,t)=\Phi+\int_s^t\Big(\A u(\cdot,r)+
\varphi(\cdot,r)\Big)\,dr+ \sum_{i=1}^N
\int_s^t(B_iu(\cdot,r)+h_i(\cdot,r))\,dw_i(r)
\label{intur} \eaa for all $t$ such that $s<t\le T$, and this
equality is satisfied as an equality in $Z_T^{-1}$.
\end{definition}
\begin{definition} 
\label{defsolltion2} \rm We say that equation (\ref{4.8}) is
satisfied for $p\in Y^1$, $\Psi\in Z_T^0$, $\chi_i\in X^0$ if
\baa 
\label{intur1} p(\cdot,t)=\Psi+\int_t^T\Big(\A^*p(\cdot,s) +
\sum_{i=1}^NB_i^*\chi_i(\cdot,s)+\xi(\cdot,s)\Big)\,ds -\sum_{i=1}^N
\int_t^T\chi_i(\cdot,s)\,dw_i(s)
\eaa for any $t\in[0,T]$.  The equality here is assumed to be an
equality in the space $Z_T^{-1}$.
\end{definition}
Note that the condition on $\p D$ is satisfied in the following
sense:   $u(\cdot,t,\o)\in H^1$ and $p(\cdot,t,\o)\in H^1$ for a.e.
\ $t,\o$. Further,  $u,p\in Y^1$, and  the value of  $u(\cdot,t)$ or
$p(\cdot,t)$  is uniquely defined in $Z_T^0$ given $t$, by the
definitions of the corresponding spaces. The integrals with $dw_i$
in (\ref{intur}),\ref{intur1}) are defined as elements of $Z_T^0$.
The integrals with $ds$ are defined as elements of $Z_T^{-1}$.
(Definitions \ref{defsolltion}-\ref{defsolltion2} require for
(\ref{4.1}) (\ref{4.8}) that these integral are equal  to  elements
of $Z_T^{0}$ in the sense of equality in $Z_T^{-1}$).
\subsubsection*{Existence theorems and
known fundamental inequalities}
\par The following Lemma combines the first and the second fundamental inequalities and related existence result
for forward SPDEs. It gives analogs of the so-called "energy
inequalities", or "the fundamental inequalities" known for
deterministic parabolic equations (Ladyzhenskaya {\it et al}
(1969)).
\begin{lemma}
\label{lemma1} Let either $k=-1$ or $k=0$. Assume that Conditions.
In addition, assume that if $k=0$, then $\b_i(x,t,\o)=0$ for $x\in
\p D$, $i=1,...,N$ and \baaa \esssup_{\o}\sup_{(x,t)\in Q} \Bigl|
\frac{\p^2 b}{\p x_k\p x_m}(x,t,\o)\Bigr|< +\infty. \eaaa  Let
$\varphi\in X^{k}(s,T)$, $h_i\in X^{k+1}(s,T)$, and $\Phi\in
Z_s^{k+1}$. Then problem (\ref{4.1}) has an unique solution $u$ in
the class $Y^1(s,T)$, and the following analog of the first
fundamental inequality is satisfied:  \be \label{4.2} \| u
\|_{Y^{k+2}(s,T)}\le c \left(\| \varphi \|
_{X^{k}(s,T)}+\|\Phi\|_{Z^{k+1}_s}+ \sum_{i=1}^N\|h_i
\|_{X^{k+1}(s,T)}\right), \ee where  $c=c({\cal P}_1)$ is a constant
that depends on ${\cal P}_1 $ only.
\end{lemma}
\par
The statement  of Lemma  \ref{lemma1}  for $k=-1$ corresponds to the
first fundamental inequality; it is a special case of Theorem 3.4.1
from Rozovskii (1990). The statement for $k=0$ corresponds to the
second energy inequality; it was obtained in Dokuchaev (2005).
\par
The following Lemma gives the first  fundamental inequalities and
related existence result for backward SPDEs.
\begin{lemma} 
\label{Th4.2} [Dokuchaev (1992,2010a)] For any $\xi \in X^{-1}$ and
$\Psi\in Z_T^0$, there exists a pair $(p,\chi)$ such that $p\in
Y^1$, $\chi=(\chi_1,\ldots, \chi_N)$, $\chi_i\in X^0$, and
(\ref{4.8}) is satisfied. This pair is uniquely defined, and the
following analog of the first fundamental inequality is satisfied:
\be\label{1FE} \|p\|_{Y^1}+\sum_{i=1}^N\|\chi_i\|_{X^0}\le
c(\|\xi\|_{X^{-1}}+\|\Psi\|_{Z_T^{0}}), \ee where $ c=c({\cal
P_1})>0$ as a constant that does not depend on $\xi$ and $\Psi$.
\end{lemma}
Therefore, only the second fundamental inequality for backward SPDEs
is missed.
\section{The main result:
the second fundamental inequality for backward equations}
\label{secMain}
\par Starting from now, we assume that the following addition
conditions are satisfied.
\begin{condition} \label{condA}  There exists a constant $\d>0$ such that
\baa
 \label{Main1'} \sum_{i=1}^Ny_i^\top  b
(x,t,\o)\,y_i-\frac{1}{2}\left(\sum_{i=1}^Ny_i^\top\b_i(x,t,\o)\right)^2
\ge \d_1\sum_{i=1}^N|y_i|^2\nonumber\\ \quad \forall\,
\{y_i\}_{i=1}^N\subset \R^n,\ (x,t)\in D\times [0,T],\ \o\in\O. \eaa
\end{condition}
\par
For an integer $M>0$, let $\Theta_b(M)$ denote the class of all
matrix functions  $b$  such that all conditions imposed in Section
\ref{SecC} are satisfied, and there exists a set
$\{t_k\}_{i=0}^M=\{t_k(M)\}_{i=0}^M$ such that $0=t_0<t_1<\cdots
<t_M=T$ such that $\max_k|t_k-t_{k-1}|\to 0$ as $M\to +\infty$, and
that the function $b(x,t,\o)=b(x,\o)$  does not depend on $t$ for
$t\in[t_i,t_{i+1})$. In particular, this means that   $b(x,t,\cdot)$
is $\F_{t_i}$-measurable for all $x\in D$, $t\in [t_i,t_{i+1})$.
\par Set
$\Theta_b\defi \cup_{M>0} \Theta_b(M)$.
\par
The following Condition \ref{condTheta} is rather technical.
\begin{condition}
\label{condTheta} The matrix $b$ is such that all conditions imposed
in Section \ref{SecC} are satisfied, and that there exits a sequence
$\{b^{(M)}\}_{M=1}^{+\infty}\subset\Theta_b$ such  at leats one of
the following conditions is satisfied:
\begin{enumerate}
\item
$\|b^{(M)}-b\|_{\W^1_{\infty}}\to 0$ as $M\to +\infty$.
\item
 Condition \ref{condA} is
satisfied for $b$ replaced by $b^{(M)}$, with the same $\d_1>0$ for
all $M$, and
$\|b^{(M)}(\cdot,t,\o)-b(\cdot,t,\o)\|_{W_{\infty}^1(D)}\to 0$ for
a.e. (almost every) $(t,\o)$ as $i\to +\infty$.
\end{enumerate}\end{condition}
\par
 We denote by $\oo\Theta_b$ the class of all functions such  $b$
that Condition \ref{condTheta} is satisfied.
\par
To proceed further, we assume that Conditions
\ref{cond3.1.A}-\ref{cond3.1.B}  remain in force starting from here
and up to the end of this paper, as well as the previously
formulated conditions.
\par
Let ${\cal P}=\{{\cal P}_1,\d_1\}$.
\begin{theorem} 
\label{Th2FE} For any $\xi \in X^0$ and $\Psi\in Z_T^1$, there
exists a pair $(p,\chi)$, such that $p\in Y^2$,
$\chi=(\chi_1,\ldots, \chi_N)$, $\chi_i\in X^1$ and (\ref{4.8}) is
satisfied. This pair is uniquely defined, and the following analog
of  the second fundamental inequality holds: \be\label{2FE}
\|p\|_{Y^2}+\sum_{i=1}^N\|\chi_i\|_{X^1}\le
c(\|\xi\|_{X^{0}}+\|\Psi\|_{Z_T^{1}}), \ee where $ c>0$ is a
constant that depends only on ${\cal P}$.
\end{theorem}
 Repeat that estimate
(\ref{2FE}) represents an analog of  the second fundamental
inequality.
\subsubsection*{On the strengthened coercivity condition}
Let us discuss the properties  Condition \ref{condA} and compare it
with Condition \ref{cond3.1.A}. First, let us note  that it can
happen that Condition \ref{cond3.1.A} holds but Condition
\ref{condA} does not hold. It can be seen from the following
example.
\begin{example}\label{exa1}{\rm
Assume that $n=2$, $N=2$,\baaa
 \b_1\equiv\Biggl(
                         \begin{array}{c}
                           1 \\
                           0 \\
                         \end{array}
                       \Biggr),\quad
 \b_2\equiv \Biggl(
                         \begin{array}{c}
                           0 \\
                           1 \\
                         \end{array}
                       \Biggr), \qquad b\equiv
\frac{1}{2}(\b_1\b_1^\top+\b_2\b_2^\top)+0.01I_2=0.51I_2, \eaaa
where $I_2$ is the unit matrix in $\R^{2\times 2}$. Obviously,
Condition \ref{cond3.1.A} holds. However, Condition \ref{condA} does
not hold for this $b$;  to see this, it suffices to take $y_1=\b_1$
and $y_2=\b_2$. }\end{example}
\begin{remark}{\rm
Assume that the estimate in Condition \ref{cond3.1.A} holds with
$\d=0$ only (i.e., the forward equation is dissipative, in the terms
of Rozovskii (1990)). This important model covers Kolmogorov type
equations for conditional densities of non-Markov Ito processes (see
Rozovsky (1990), Chapter 6, and Dokuchaev (1995)). If we approximate
the operator $\A$ by the operator $\A+\e\Delta$, where $\Delta$ is
the Laplacian, then Condition \ref{cond3.1.A} holds for the new
operator for arbitrarily small $\e>0$. This approximation of a
dissipative equation by a coercive one  is a useful tool for
investigation of dissipative equations and distributions of
non-Markov Ito processes. Example \ref{exa1} shows that,
unfortunately, general dissipative equations cannot be approximated
by equations such that  Condition \ref{condA} holds. }\end{remark}
\par
The following theorems clarify the relations between Conditions
\ref{condA} and \ref{cond3.1.A}.
\begin{theorem}\label{propNewOld}
If Condition \ref{condA} holds then Condition \ref{cond3.1.A} holds.
\end{theorem}
Let us give some useful criterions of validity of Condition
\ref{condA}.
\begin{theorem}\label{propON1} If
$n=1$ and  Condition \ref{cond3.1.A} holds, then
 Condition \ref{condA} holds.
\end{theorem}\begin{theorem} \label{propON2}
 Condition \ref{condA} holds
 if
 there exist $N_0\in\{1,...,N\}$  and  $\d_2>0$ such that
 $\b_i\equiv 0$ for $i>N_0$ and
\be
 \label{Main1''} y^\top  b
(x,t,\o)\,y-\frac{N_0}{2}|y^\top\b_i(x,t,\o)|^2 \ge \d_2|y|^2
\quad\forall\, y\in \R^n,\ (x,t)\in  D\times [0,T],\ \o\in\O,\
i=1,...,N_0. \ee
\end{theorem}
\begin{corollary} If $N=1$ and Condition \ref{cond3.1.A} holds
then Condition \ref{condA} holds.
\end{corollary}
\section{Some
applications}\label{Secappl} So far, the main application is the
representation theorem  for functionals of non-Markov processes and
their first exit times from bounded domains. These functionals are
represented via solutions of backward parabolic Ito equations. The
previously known results about regularity of the solution of the
backward SPDE for $p$ were insufficient for the case of domains with
boundary, and the representation result was never before obtained
for this case. It was done only using the additional regularity in
the form of the second fundamental inequality given in Theorem
\ref{Th2FE} (Dokuchaev (2010b)). Therefore, this regularity result
opens ways to systematics of first exit times of non-Markov
processes.
\par
In addition, a priori estimates obtained above helps to establish w
can show that the solution of (\ref{4.8}) is robust with respect to
small in $L_{\infty}$ norm disturbances of the coefficients.
\par
 Consider two problems (\ref{4.8}), with coefficients
 $$(b,f,\lambda,\xi,\b_i,\oo\b_i,\Psi)
 =(b^{(k)},f^{(k)},\lambda^{(k)},\xi^{(k)},\b_i^{(k)},\oo\b_i^{(k)},\Psi^{(k)}), \quad k=1,2,$$
 such that Conditions  \ref{cond3.1.A}-\ref{cond3.1.B}  and   \ref{condA}-\ref{condTheta}
 are satisfied for both sets of functions. Let ${\cal P}^{(k)}$ be the corresponding sets of
 parameters.  Let $(p^{(k)},\chi^{(k)}_1,...,\chi^{(k)}_N  )$ be the corresponding solutions of problem
(\ref{4.8}), $k=1,2$.
 \begin{theorem}\label{ThRobust}
There exists a constant $c=c({\cal P}^{(1)},\|u^{(2)}\|_{Y_2})$ such
that   $$\| p^{(1)}-p^{(2)}
\|_{{Y}^2}+\sum_{i=1}^N\|\chi^{(1)}_i-\chi^{(2)}_i\|_{X^1}
 \le   c M,$$ where
 \baaa
 &&M\defi\esssup_{x,t\o}\biggl( |b^{(1)}(x,t,\o)-
b^{(2)}(x,t,\o)|
 +|f^{(1)}(x,t,\o)-f^{(2)}(x,t,\o)|\\&&
+|\lambda^{(1)}(x,t,\o)-\lambda^{(2)}(x,t,\o)|+ \sum_{i=1}^N|
\b_i^{(1)}(x,t,\o)- \b_i^{(2)}(x,t,\o)|
 \\&& + \sum_{i=1}^N| \oo\b_i^{(1)}(x,t,\o)- \oo\b_i^{(2)}(x,t,\o)|
\biggr) +\|\xi^{(1)}-\xi^{(2)}\|_{X^0} +\|\Psi^{(1)}-\Psi^{(2)}
\|_{Z_0^1}  .\nonumber\label{rob1}
 \eaaa
\end{theorem}
\par
Note that the first fundamental inequality  can help to establish
robustness  only with respect to deviations of $b$ that are small
together with  their derivatives in $x$, and this restriction is
necessary even for robustness in $X^0$. Theorem \ref{ThRobust}
establishes robustness in $Y^2$ for the disturbances of the
coefficients that are small in $L_{\infty}$-norm only. For instance,
if $b$ is replaced for $b+\xi$, where
$\esssup_{x,t,\o}|\xi(x,t,\o)|\le \e$ for a small $\e>0$, then
Theorem \ref{ThRobust} ensures that the corresponding solution of
(\ref{4.8}) is close  in  $Y^2$ to the original one. \vspace{4mm}
 \par
 The rest part of the
paper is devoted to the proofs of results given above.
\section{Auxiliary facts for backward equations}
In this section, we collect some facts that will be used for the
proof of Theorem \ref{Th2FE}. Lemmas
\ref{Th4.2operators}-\ref{lemmaMark} given below were obtained in
Dokuchaev (2010a), where the their proof can be found.
\subsection{Decomposition of operators $L$ and
$\M_i$}\label{SecDecomposition}
 Introduce  operators $L(s,T):X^{-1}(s,T)\to
Y^1(s,T)$, $\M_{i}(s,T):X^{0}(s,T)\to Y^1(s,T)$, and
$\L(s,T):Z^0_s\to Y^1(s,T)$, such that
$$u=L(s,T)\varphi+\L(s,T)\Phi+\sum_{i=1}^N\M_{i}(s,T)h_i,$$ where
$u$ is the solution in $Y^1(s,T)$ of   problem (\ref{4.1}). These
operators are linear and continuous; it follows immediately from
Lemma \ref{lemma1}. We will denote by $L$, $\M_i$, and $\L$, the
operators $L(0,T)$, $\M_i(0,T)$, and $\L(0,T)$, correspondingly.
\par
For $t\in[0,T]$, define operators $\I_t: C([0,T];Z_T^{k})\to
Z^{k}_t$ such that $\I_tu=u(\cdot,t)$.
\begin{lemma}\label{Th4.2operators}
In the notations of Lemma \ref{Th4.2}, the following duality
equation is satisfied: \baa p={\LU}^*\xi+(\I_T\LU)^*\Psi,\quad
\chi_i=\M_i^*\xi+(\I_T\M_i)^*\Psi,\quad
p(\cdot,0)=\L^*\xi+(\I_T\L)^*\Psi,\label{duality}\eaa where
${\LU}^*: X^{-1}\to X^1$, ${\M}_i^*: X^{-1}\to X^0$, $(\I_T\LU)^*:
Z_0^{0}\to X^1$, $(\I_T\M_i)^*: Z_0^{0}\to X^0$, and $(\I_T\L)^*:
Z_T^{0}\to Z_0^0$,
 are the operators that are adjoint  to the operators ${\LU}:
X^{-1}\to X^1$, $\M_i : X^{0}\to X^1$, $\I_T\M_i: X^{-1}\to Z_T^0$,
$\I_T\M_i: X^0\to Z_T^{0}$, and $\,\I_T\L: Z_0^0\to Z_T^{0}$,
respectively.
\end{lemma}
\par
Our method of  proof of fundamental inequalities is based on
decomposition of the operators to superpositions of simpler
operators.
\begin{definition} 
\label{Def4.2} \rm Define operators $\KK: Z_0^{0}\to Y^1$, $\Q_0:
X^{-1}\to Y^1$, $\Q_i:\allowbreak X^0\to Y^1$, $i=1,...,N$, as the
operators $\L: Z_0^{0}\to Y^1$, $\LU: X^{-1}\to Y^1$,
$\M_i:\allowbreak X^0\to Y^1$, $i=1,...,N$, considered for the case
when $B_i=0$ for all $i$.
\end{definition}
By Lemma \ref{lemma1}, these linear operators are continuous. It
follows from the definitions  that
$$
\KK\Phi+\Q_0\eta+\sum_{i=1}^N\Q_ih_i=V,
$$
where $\eta\in X^{-1}$, $\Phi\in Z_0^0$, and $h_i\in X^0$, and where
$V$ is the solution of the problem
\begin{eqnarray}
\begin{array}{c}
\label{4.11} d_tV=(\A V+\eta)\,dt + \sum_{i=1}^N h_i\,dw_i(t),
\\
 V|_{t=0}=\Phi,\quad V(x,t,\o)\,|_{x\in \p D}=0.
\end{array}
\end{eqnarray}
 Define the operators \be\label{P*} P\defi\sum_{i=1}^N\Q_iB_i,
\quad P^*\defi\sum_{i=1}^NB_i^*\Q_i^*. \ee By the definitions, the
operator $P: X^{1}\to X^{1}$ is continuous, and $P^*: X^{-1}\to
X^{-1}$ is its adjoint  operator. Hence the operator $P^*: X^{-1}\to
X^{-1}$ is continuous.
 Let
$$ P_0\defi \I_T\sum_{i=1}^N\Q_iB_i,\quad P_0^*\defi
\sum_{i=1}^NB_i^*(\I_T\Q_i)^*. $$  By the definitions, the operator
$P_0: X^{1}\to Z_T^{0}$ is continuous, and $P_0^*: Z_T^{0}\to
X^{-1}$ is its adjoint  operator. Hence the operator $P_0^*: Z_T^\to
X^{-1}$ is continuous.
\begin{lemma} 
\label{prop4.3} The operator $(I-P)^{-1}: X^{1 }\to X^{1}$ is
continuous, and \baa &{\LU}=(I-P)^{-1}\Q_0,
 &\M_i=(I-P)^{-1}\Q_i,\nonumber\\
&\I_T\LU=P_0(I-P)^{-1}\Q_0+\I_T\Q_0,
&\I_T\M_i=P_0(I-P)^{-1}\Q_i+\I_T\Q_i, \label{I-P}\eaa $i=1,...,N$.
The operator $(I-P^*)^{-1}: X^{-1}\to X^{-1}$ is also continuous,
and \baa &\LU^*=\Q_0^*(I-P^*)^{-1},\quad
&\M_i^*=\Q_i^*(I-P^*)^{-1},\nonumber
\\
&(\I_T\LU)^*=\Q_0^*(I-P^*)^{-1}P_0^*+(\I_T\Q_0)^*, &
(\I_T\M_i)^*=\Q_i^*(I-P^*)^{-1}P_0^*+(\I_T\Q_i)^*.\hphantom{xxxx}
\label{I-P*} \eaa \end{lemma} In fact, Lemma \ref{prop4.3} allows to
split  represent solution (\ref{4.8}) via solution of much simpler
problem with $B_i\equiv 0$ and via inverse operator $(I-P^*)^{-1}$.
It can be illustrated as the following.
\begin{corollary}\label{corrgxi}\begin{enumerate}
\item For $\Psi=0$, the solution $(p,\chi_1,...,\chi_N)$ of problem
(\ref{4.8}) can be represented as $p=\Q_0^*g$, $\chi_i=\Q_i^*g$,
where $g=\xi+ \sum_{i=1}^NB_i^*\chi_i$, and where
$\sum_{i=1}^NB_i\chi_i=P^*g$.
\item For general $\Psi$, the solution $(p,\chi_1,...,\chi_N)$ of problem
(\ref{4.8}) can be represented as \baaa p=\Q_0^*g+(\I_T\Q_0)^*\Psi,
\quad \chi_i=\Q_i^*g+(\I_T\Q_i)^*\Psi,\eaaa
 where $g=\xi+
\sum_{i=1}^NB_i^*\chi_i$, and where
$\sum_{i=1}^NB_i\chi_i=P^*g+P_0^*\Psi$. In other words,
$g=(I-P^*)^{-1}\xi+ (I-P^*)^{-1}P_0^*\Psi.$
\end{enumerate}
\end{corollary}
It appears that this representation helps to establish the second
fundamental inequality.
\subsection{Semi-group property for backward equations} It is known that the
forward SPDE is casual (or it has semigroup property): if
$u=L\varphi+\L\Phi$, where $\varphi\in X^{-1}$, $\Phi\in Z_0^0$,
then \be\label{casF}
u|_{t\in[\t,s]}=L(\t,s)\varphi+\L(\t,s)u(\cdot,\t).\ee To proceed
further, we need  a similar property for the backward equations.
\begin{lemma}\label{lemmaMark}
Let $0\le \t< s<T$, and let $p=L^*\xi$, $\chi_i=\M_i\xi$, where
$\xi\in X^{-1}$ and $\Psi\in Z_T^0$. Then \baa
&&p|_{t\in[\t,s]}=L(\t,s)^*\xi|_{t\in[\t,s]}+(\I_sL(\t,s))^*p(\cdot,s),\label{p-cas}
\\
&&p(\cdot,\t)=(\I_{\t}\L(\t,s))^*p(\cdot,s)+\L(\t,s)^*\xi,\label{p-cas0}
\\
&&\chi_k|_{t\in[\t,s]}=\M_k(\t,s)^*\xi|_{t\in[\t,s]}+(\I_s\M_i(\t,s))^*p(\cdot,s),
\quad k=1,...,N. \label{chi-cas} \eaa
\end{lemma}
\par
Note that this semi-group property implies causality for backward
equation (which is a non-trivial fact due the presence of $\chi$).
 \subsection{A special estimate for deterministic PDEs}
We use notations $ \nabla u\defi \Bigl(\frac{\p u}{\p x_1},\frac{\p
u}{\p x_2},...,\frac{\p u}{\p x_n}\Bigr)^\top, $ for functions
$u:\R^n\to\R$. In addition, we use the notation $ (u,v)_{H^0}\defi
\sum_{i=1}^n(v_i,u_i)_{H^0}$ for functions $u,v: D\to\R^n$, where
$u=(u_1,...,u_n)$ and $v=(v_1,...,v_n)$.
\par
For $u\in H^1$, let \baa \|u\|_{\w H^1(t,\o)}\defi (\nabla
u,b(\cdot,t,\o)\nabla u)^{1/2}= \Bigl(\sum_{i,j=1}^n \int_D \frac{\p
u}{\p x_i }(x)b_{ij}(x,t,\o)\frac{\p u}{\p x_j}(x)dx
\Bigr)^{1/2}.\label{wH} \eaa
\par
For $K>0$, introduce the operator $\A^*_K=\A^*-KI$, i.e., $\A^*_K
u=\A^* u-Ku$.
\begin{lemma}\label{lemma<<1} Let $\t,\tau\in
[0,T]$ be given,  $0\le\t<\tau\le T$. Let the function
$b(x,t,\o)=b(x,\o)$ be constant in $t\in[\t,\tau]$ for a.e. $x$,
$\o$.  Let $h=h(x,t,\o)\in L_2(D\times[\t,\tau])$, and let
$u=u(x,t,\o):D\times[\t,\tau]\times\O\to\R$ be the solution of the
boundary value problem
 \be\ba \frac{\p u}{\p
t}+\A_K^*u=-h,\quad t\in(\t,\tau)\\
u(x,\tau) =0,\quad u(x,t)|_{x\in \p D}=0, \label{oog3}\ea\ee
 Then
for any $\e>0$, $M>0$, there exists $K=K(\e,M,{\cal P})>0$ such that
\baaa \sup_{t\in[\t,\tau]}\|u(\cdot,t,\o)\|_{\ww H^1(t,\o)}^2+
M\sup_{t\in[\t,\tau]}\|u(\cdot,t,\o)\|^2_{H^0}\le
\frac{1+\e}{2}\int_\t^\tau\|h(\cdot,t,\o)\|_{H^0}^2dt\quad\hbox{a.s.}
\eaaa
\end{lemma}
This lemma follows immediately from Theorem 1 and Corollary 1 from
Dokuchaev (2008).
\section{The
proof of Theorem \ref{Th2FE}} By Lemma \ref{Th4.2operators}, it
suffices to show that the operators $L^*: X^0\to Y^2$,
$(\I_TL)^*:Z_T^1\to Y^2$, and $\M_i^*:X^0\to X^1$,
$(\I_T\M_i)^*:Z_T^1\to X^1$, are continuous, and that their norms
are less or equal than a constant $c=c({\cal P})$.
\par
We define the operators $L^*(s,T)$,  $\M_i^*(s,T)$, $
(\I_TL(s,T))^*$, and $(\I_T\M_i(s,T))^*$, similarly to $L^*$,
$\M_i^*$, $(\I_TL)^*$, and $(\I_T\M_i)^*$, with time interval
$[0,T]$ replaced by $[s,T]$.
\par
We denote by $\oo{{\cal P}_T}$ the completion (with respect to the
measure $\oo\ell_1\times\P$ of the $\s$-algebra of subsets of
$[0,T]\times\O$, generated by functions that are progressively
measurable with respect to $ \oo{{\cal B}}_{1}\times\F_T$. Let $ \oo
X^{k}\defi L^{2}\bigl([0,T ]\times\Omega, {\oo{\cal P
}}_T,\oo\ell_{1}\times\P;  H^{k}\bigr). $
\par
 Let ${\cal E}$ be the
operator of projection of $\oo X^1$ onto $X^1$.
\par
Let $\xi\in X^0$, $\Psi\in Z_T^1$,  and let $\oo p$ be the solution
of the boundary value
problem in $Q$  \be\ba \frac{\p \oo p}{\p t}+\A^*\oo p=-\xi,\quad t\le T,\\
\oo p|_{t=T}=\Psi,\quad \oo p(x,t,\o)|_{x\in \p D}=0.
\label{oop2}\ea\ee  By the second fundamental inequality for
deterministic parabolic equations, it follows that the solution of
(\ref{oop2}) is such that $\oo p\in \oo X^2\cap \CC^1$, (\ref{2fee})
holds and \baa \|\oo p\|_{\oo X^2}+\|\oo p\|_{\CC^1} \le
c\left(\|\xi\|_{X^{0}}+ \|\Psi\|_{Z_T^{0}}\right), \label{2fee}\eaa
where $c=c({\cal P})>0$ is a constant.   This fact is well known; if
the function $b(x,t,\o)$ is almost surely continuous, then
(\ref{2fee}) follows Theorem IV.9.1 from Ladyzenskaya {\em et al}
(1968).  Since the derivative $\p b/\p x$ is bounded, the condition
that $b$ is continuous can be lifted. In this case, (\ref{2fee})
follows from Theorem 3.1 from Dokuchaev (2005).
\par
 By Martingale Representation Theorem, there exist
functions $\g_i(\cdot,t,\cdot)\in X^0$ such that \baa \oo
p(x,t,\o)=\E\{\oo
p(x,t,\o)|\F_0\}+\sum_{i=1}^N\int_0^T\g_i(x,t,s,\o)dw_i(s).\label{clarkp}
\eaa
\begin{lemma}\label{lemmaQ}  Assume that the function  $\mu=(b, f,\lambda)$ is such that $\mu(x,t,\o)$ is
$\F_0$-measurable for all $x\in D$.    Let  $\xi\in X^0$, $\Psi\in
Z_T^1$,   let $\oo p$ be the solution of (\ref{oop2}), and let
$\g_j$ be the processes presented in (\ref{clarkp}). Let
$p,\chi_1,...,\chi_2$ be defined as  \baa p\defi{\cal E}\oo p,\quad
\chi_i(x,s,\o)\defi \g_i(x,s,s,\o).\label{pchi}\eaa Then $p\in Y^1$,
$\chi_i\in X^1$, and \baa\|p\|_{Y^2}+\sum_{i=1}^N\|\chi_i\|_{X^1}
\le c\left(\|\xi\|_{X^{0}}+ \|\Psi\|_{Z_T^{0}}\right),
\label{pest}\eaa where $c=c({\cal P})>0$ is a constant. In addition,
\be \label{chi2} p=\Q_0^*\xi+(\I_T\Q_0)^*\Psi,\quad
\chi_i=\Q_i^*\xi+(\I_T\Q_i)^*\Psi. \ee
\end{lemma}
\par
{\em Proof of Lemma \ref{lemmaQ}.} By Martingale Representation
Theorem, there exist functions $\g_i(\cdot,t,\cdot)\in X^0$,
$\g_{\xi i}(\cdot,t,\cdot)\in X^0$, and $\g_{\Psi i}(\cdot)\in X^1$,
such that (\ref{clarkp}) holds as well as  \baaa
&&\xi(x,t,\o)=\E\{\xi(x,t,\o)|\F_0\}+\sum_{i=1}^N\int_0^T\g_{\xi
i}(x,t,s,\o)dw_i(s),  \\
&&\Psi(x,\o)=\E\{\Psi(x,\o)|\F_0\}+\sum_{i=1}^N\int_0^T\g_{\Psi
i}(x,s,\o)dw_i(s). \eaaa Moreover, it follows that ${\cal D}
g_i(\cdot,t,\cdot)\in X^0$, where either ${\cal D}\g =\p\g/\p t$ or
${\cal D} \g=\A^*\g$, and \baaa {\cal D}\oo p(x,t,\o)=\E\{{\cal
D}\oo p(x,t,\o)|\F_0\}+\sum_{i=1}^N\int_0^T{\cal
D}\g_i(x,t,s,\o)dw_i(s). \eaaa
\par
By (\ref{oop2}), it follows that  \be\ba \frac{\p \g_i}{\p
t}(\cdot,t,s,\o)+\A^*\g_i(\cdot,t,s,\o)=-\g_{\xi
i}(\cdot,t,s,\o),\quad t\in (0,T),
\\ \g_i(x,T,s,\o) =\g_{\Psi i}(x,s,\o),\qquad \g_i(x,t,s,\o)|_{x\in \p
D}=0. \label{oog2}\ea\ee
\par
Again, it follows from the second fundamental inequality for
deterministic parabolic equations that \baaa \sup_{t\in[s,T]}\|\g_{
i}(\cdot,t,s,\o)\|_{H^1}^2\le c\left(\int_s^T\|\g_{\xi
i}(\cdot,t,s,\o)\|_{H^{0}}^2dt+\|\g_{\Psi
i}(\cdot,s,\o)\|_{H^1}^2\right), \eaaa where $c=c(T,n,D)>0$ is a
constant. Hence \baaa \|\g_{ i}(\cdot,s,s,\o)\|_{H^1}^2\le
c\left(\int_s^T\|\g_{\xi i}(\cdot,t,s,\o)\|_{H^{0}}^2dt+\|\g_{\Psi
i}(\cdot,s,\o)\|_{H^1}^2\right).
 \eaaa  This estimate together with (\ref{2fee}) ensures that
(\ref{pest}) holds for $p$ and $\chi_i$ defined by (\ref{pchi}).
\par
Let us show that (\ref{chi2}) holds.
\par
 Clearly, \baaa \oo p(x,t,\o)=
p(x,t,\o)+\sum_{i=1}^N\int_t^T\g_i(x,t,s,\o)dw_i(s), \eaaa and \baaa
\oo p(\cdot,t)=\int_t^T\Bigl(\A^*\oo p(\cdot,s)
+\xi(\cdot,s)\Bigr)\,ds.\eaaa Hence \baaa
p(\cdot,t)&=&\Psi+\int_t^T\Bigl(\A^*p(\cdot,s)
+\xi(\cdot,s)\Bigr)\,ds \\
&&+\sum_{i=1}^N\biggl[\int_t^Tds\int_s^T[
\A^*\g_i(\cdot,s,r)+\g_{\xi i}(\cdot,s,r)]dw_i(r)-\int_t^T
\g_i(\cdot,t,s)dw_i(s)\biggr]
\\
&&=\Psi+\int_t^T\Bigl(\A^*p(\cdot,s) +\xi(\cdot,s)\Bigr)\,ds \\
&&+\sum_{i=1}^N\biggl[\int_t^T dw_i(r)\int_t^r[
\A^*\g_i(\cdot,s,r)+\g_{\xi i}(\cdot,s,r)]ds-\int_t^T \g_i(\cdot,t,s)dw_i(s)\biggr]\\
\\
&&=\Psi+\int_t^T\Bigl(\A^*p(\cdot,s)+
\xi(\cdot,s)\Bigr)\,ds\\&&+\sum_{i=1}^N\int_t^T
dw_i(s)\biggl[\int_t^s[ \A^*\g_i(\cdot,r,s)+\g_{\xi
i}(\cdot,r,s)]dr-\g_i(\cdot,t,s)\biggr].\eaaa By (\ref{oog2}),
$$
\g_i(\cdot,t,s)-\int_t^s[ \A^*\g_i(\cdot,r,s)+\g_{\xi
i}(\cdot,r,s)]dr =\g_i(\cdot,s,s).
$$
By  (\ref{pchi}), we have selected
$\g_i(\cdot,s,s)=\chi_i(\cdot,s)$. It follows that \baaa
p(\cdot,t)=\Psi+\int_t^T\Biggl(\A^*p(\cdot,s)+\xi(\cdot,s)\Biggr)\,ds
-\,\sum_{i=1}^N \int_t^T\chi_i(\cdot,s)\,dw_i(s). \eaaa Finally, we
obtain (\ref{chi2}) from Lemma \ref{Th4.2operators} applied to the
operators $\Q_0^*$, $\Q_i^*$, $(\I_T\Q_0)^*$, and $(\I_T\Q_i)^*$,
$i=1,...,N$, considered as special cases of $\LU^*$, $\M_i^*$,
$(\I_T L_0)^*$, and $(\I_T\M_i)^*$, respectively. This completes the
proof of Lemma
 \ref{lemmaQ}. $\Box$
\vspace{4mm}
 \par
 In the following proof, we will  explore the following observation:
  if $\lambda$ is replaced  by $
\lambda^{(K)}(x,t,\o)\defi \lambda(x,t,\o)+K$,  i.e., if $\A$ is
replaced by $\A_K=\A v+ KI$, then the solution $u$ of the problem
(\ref{4.1}) has to be replaced by the process $$u(x,t,\o)e^{-Kt},$$
and
 the solution $(p,\chi_1,...,\chi_N)$
of the problem (\ref{4.8}) has to be replaced by the process
$$\Big(p(x,t,\o)e^{K(T-t)},\chi_1(x,t,\o)e^{K(T-t)},...,\chi_N(x,t,\o)e^{K(T-t)}\Big).$$
 Therefore, it suffices to prove theorem for any case when
$\lambda$ is replaced  for $\lambda^{(K)}(x,t,\o)\defi
\lambda(x,t,\o)+K$ with some  $K>0$, and this $K$ can be taken
arbitrarily large.
\par
For linear normed spaces $\X$ and $\Y$, we denote by $\|{\cal T}
\|_{\X,\Y}$ the norm of an operator ${\cal T}: \X\to \Y$.
\begin{lemma}\label{lemma<<<1} Let $0\le s<T$, and let the function
$\mu=(b,f,\lambda)$ be such that $\mu(x,t,\cdot)$ is
$\F_{s}$-measurable for all $x\in D$, $t\in [s,T)$. Moreover, we
assume that $b(x,t,\o)=b(x,\o)$ does not depend on $t\in[s,T]$.
 Then there exist $K>0$ such that if $\lambda$ is replaced  by
$\lambda(x,t,\o)+K$, then \baaa \|L^*(s,T)\|_{X^0(s,T),
Y^2(s,T)}+\|(\I_TL(s,T))^*\|_{Z_T^1,
Y^2(s,T)}+\sum_{i=1}^n\|\M_i^*(s,T)\|_{X^0(s,T), X^1}\\
+\sum_{i=1}^n\|(\I_T\M_i(s,T))^*\|_{Z_T^1, X^1(s,T)} \le c, \eaaa
where $c\in(0,+\infty)$ depends only on $K$ and ${\cal P}$.
\end{lemma}
\par
{\it Proof of Lemma \ref{lemma<<<1}.} To simplify the notations, we
consider only the case when $s=0$.
\par
By (\ref{duality}) and (\ref{I-P*}), it suffices to show that the
operator $(I-P^*)^{-1}:X^0\to X^0$ is continuous.  For this, it
suffices to show that there exist $K>0$ such that if $\lambda$ is
replaced for $\lambda(x,t,\o)+K$ then $\|P^*\|_{X^0, X^0}<1$.
\par
Let  $\xi\in X^0$, let $\oo p$ be the solution of (\ref{oop2}), and
let $\g_j$ be the processes presented in (\ref{clarkp}) with
$\Psi=0$. Let $p,\chi_1,...,\chi_2$ be defined by (\ref{pchi}) with
$\Psi=0$. In this case, \be\ba \frac{\p \g_i}{\p
t}(\cdot,t,s,\o)+\A^*\g_i(\cdot,t,s,\o)=-\g_{\xi
i}(\cdot,t,s,\o),\quad t\in (0,T),
\\ \g_i(x,T,s,\o) =0,\qquad \g_i(x,t,s,\o)|_{x\in \p
D}=0. \label{oog2<<}\ea\ee
\par
By Lemma \ref{lemma<<1} applied to boundary value problem
(\ref{oog2<<}), for any $\e>0$, $M>0$, there exists $K=K(\e,M,{\cal
P})>0$ such that   \baaa \sup_{t\in[s,T]}
\|\g_i(\cdot,t,s,\o)\|_{\ww H^1(t,\o)}^2+M
\sup_{t\in[s,T]}\|\g_i(\cdot,t,s,\o)\|_{H^0}^2\le
\frac{1+\e}{2}\int_s^T\|\g_{\xi i}(\cdot,t,s,\o)\|_{H^{0}}^2dt
\quad\hbox{a.s.} \eaaa Here $\|\cdot\|_{\ww H^1(t,\o)}$ is defined
by (\ref{wH}). Hence \baaa \int_0^T\!\|\g_i(\cdot,s,s,\o)\|_{\ww
H^1(t,\o)}^2ds+M \int_0^T\!\|\g_i(\cdot,s,s,\o)\|_{H^0}^2ds\le
\frac{1+\e}{2}\int_0^T\! ds\!\int_s^T\|\g_{\xi
i}(\cdot,t,s,\o)\|_{H^{0}}^2dt. \eaaa Note that
 $$
\E\sum_{i=1}^N\int_0^Tdt\int_0^T\|\g_{\xi
i}(\cdot,t,s,\o)\|_{H^{0}}^2ds\le \|\xi \|_{X^{0}}^2.
 $$
Hence  \baaa \E\int_0^T\|\g_i(\cdot,s,s,\o)\|_{\ww H^1(t,\o)}^2ds+M
\E\int_0^T\|\g_i(\cdot,s,s,\o)\|_{H^0}^2ds\le \frac{1+\e}{2} \|\xi
\|_{X^{0}}^2 . \eaaa
 By (\ref{chi2}), it can be
rewritten as   \baa
 \E\int_0^T\|\chi_i(\cdot,t,\o)\|_{\ww H^1(t,\o)}^2dt+M
\E\int_0^T\|\chi_i(\cdot,t,\o)\|_{H^0}^2dt \le  \frac{1+\e}{2}\|\xi
\|_{X^{0}}^2 . \label{chi<}\eaa
\par
Remind that\baaa P^*\xi=\sum_{j=1}^N B_j^*\chi_j.\eaaa
\par
By Condition {\ref{condA}, there exists $M=M({\cal P})>0$ such that
\be\label{M} \biggl\|\sum_{j=1}^N B_j^*\chi_j\biggr\|^2_{H^0}\le 2
\sum_{j=1}^N\|\chi_j\|_{\ww H^1(t,\o)}^2+2M
\sum_{j=1}^N\|\chi_j\|_{H^0}^2-2\d_1\sum_{j=1}^N\|\nabla
\chi_j\|^2_{H^0}\quad \forall t,\o. \ee \par By (\ref{chi<}) and
(\ref{M}), it follows that a small enough $\e>0$ and a large enough
$K>0$ can be found such that $$ \left\|P^*\xi\right\|^2_{X^0}
=\biggl\|\sum_{i=1}^N B_i^*\chi_i\biggr\|^2_{X^0}\le
c\|\xi\|^2_{X^{0}}
$$ for this $K$ with some $c=c({\cal P},K)<1$. Hence $$
\|P^*\xi\|_{X^0}\le \sqrt{c}\|\xi\|_{X^{0}}.
$$
  Therefore,  we have proved that there
exist $K=K({\cal P})>0$ such that if $\lambda$ is replaced for
$\lambda(x,t,\o)+K$ then $\|P^*\|_{X^0, X^0}<1$, and, therefore, the
operator $(I-P^*)^{-1}:X^0\to X^0$ is continuous. By the first
equation in (\ref{chi2}), it follows that the operator
$\Q_0^*:X^0\to Y^2$ is continuous. In addition, it follows from
(\ref{chi2}) and (\ref{chi<}) that the operators $\Q_i^*:X^0\to X^1$
are continuous.
 Then the
proof of Lemma \ref{lemma<<<1} for the special case of $\Psi=0$
follows from the  first equations for adjoint operators  in
(\ref{I-P*}).
\par
To complete the proof of Lemma \ref{lemma<<<1} for general $\Psi$,
By  (\ref{pest}), (\ref{chi2}), it follows that  it suffices to show
that the operators  $(\I_T \Q_0)^*:Z_T^1\to Y^2$ and $(\I_T
\Q_i)^*:Z_T^1\to X^1$ are continuous, $i=1,...,N$. In addition, the
upper bound of the norms of these operators depends on ${\cal P}$
only. Then the proof follows from the last two equations for the
adjoint operators in (\ref{I-P*}). This completes the proof of Lemma
 \ref{lemma<<<1}. $\Box$
\vspace{4mm}
 \par
For an integer $M>0$, we denote by $\Theta(M)$ the class of all
functions  $\mu=(b,f,\lambda)$  such that all conditions imposed in
Section \ref{SecC} are satisfied, and that there exists and a set
$\{t_i\}_{i=0}^M$ such that $0=t_0<t_1<\cdots <t_M=T$ and that the
function
$\mu(x,t,\cdot)=(b(x,t,\cdot),f(x,t,\cdot),\lambda(x,t,\cdot))$ is
$\F_{t_i}$-measurable for all $x\in D$, $t\in [t_i,t_{i+1})$ and
that the function $b(x,t,\o)=b(x,\o)$  does not depend on $t$ for
$t\in[t_i,t_{i+1})$.
\par
Let $\Theta\defi \cup_{M>0} \Theta(M)$.
\begin{lemma}\label{lemma<1} Let $(b,
f,\lambda)\in \Theta(M)$ for some $M>0$.  Then there exists $K>0$
such that if $\lambda$ is replaced  by $\lambda(x,t,\o)+K$, then $$
\|L^*\|_{X^0, Y^1}+\sum_{i=1}^n\|\M_i^*\|_{X^0,
X^1}+\|(\I_TL)^*\|_{Z_T^1, Y^1}+\sum_{i=1}^n\|(\I_T\M_i)^*\|_{Z_T^1,
X^1} \le c, $$ where $c\in(0,+\infty)$ does not depend on $M$ and
depends only on $K$ and ${\cal P}$.
\end{lemma}
\par
{\it Proof} of this lemma follows immediately from Lemma
\ref{lemmaMark} and from Lemma \ref{lemma<<<1} applied consequently
for all time intervals from the definition of $\Theta(M)$ backward
from terminal time.
\begin{corollary}\label{lemmac}  Under
assumption of Lemma \ref{lemma<1},  Theorem  \ref{Th2FE} holds and
there exists $K>0$ such that the operators  $L^*:X^0\to Y^2$,
$\M_j^*:X^0\to X^1$, $(\I_TL)^*: Z_T^1\to Y^2$, and
$(\I_T\M_i)^*:Z_T^1\to X^1$, $j=1,...,N$, are continuous, and their
norms do not depend on $M$.
\end{corollary}
\par
Up to the end of this section, we assume that $\lambda$ is replaced
for $\lambda(x,t,\o)+K$ such that the conclusion of Lemma
\ref{lemma<1} holds.
\par
Now we are in position to prove Theorem  \ref{Th2FE} for the case of
$(b,f,\lambda)$ of the general kind.
\par
 Let $M=1,2,...$, $M\to +\infty$.  Let $\e\defi M^{-1}$. By
 Condition \ref{condTheta}, there exist a subsequence of $M$ such
 that there exists
 $b_{\e}\in\Theta_b(M)$ for any $M$ with the corresponding sets
 $\{t_k\}=\{t_k(M)\}$, $0=t_0<...< t_k<t_M=T$
such that $\max_k|t_k-t_{k-1}|\to 0$ as $M\to +\infty$,
$b_{\e}(x,t,\o)=b(t_k,t,\o), \quad t\in [t_k,t_{k+1})$, and that
there exists $q,r\in[1,+\infty]$ such that
$$
b_{\e}\to b \quad \hbox{in}\quad \W^1_{q,r}\quad\hbox{as}\quad \e\to
0.$$
 Further,
we introduce functions $f_{\e}$,  $\lambda_{\e}$,
 such that
\baaa f_{\e}(x,t,\o)=\E\{f(x,t,\o)|\F_{t_k}\},\quad
\lambda_{\e}(x,t,\o)=\E\{\lambda(x,t,\o)|\F_{t_k}\}, \quad t\in
[t_k,t_{k+1}). \eaaa
\begin{proposition}\label{proptech}
Let us show that Condition \ref{condA} implies that: \begin{itemize}
\item[(a)] Condition \ref{condA} is satisfied for $b$ replaced by
$b_{\e}$, with the same $\d_1>0$ for all $\e$, and \item[(b)]
Without a loss of generality, we can assume that
$\sup_{\e>0}\|b_{\e}\|_{\W_{\infty}^1}<+\infty$.
\end{itemize}
\end{proposition}
\par
{\it Proof of Proposition \ref{proptech}}. It suffices to show that
Condition \ref{condA}(i) implies (a) and that Condition
\ref{condA}(ii) implies (b).
\par
Let us show that Condition \ref{condA}(i) implies (a). Let ${\rm
A}={\rm  A}(x,t,\o)\in\R^{nN\times nN}$ be the symmetric  matrix
that defines the quadratic form on the vectors
$Y=(y_1,...,y_N)\in\R^{nN}$ in (\ref{Main1'}), and let ${\rm  A}_\e$
be the similar matrix defined for $b=b_{\e}$. By Condition
\ref{condA}, the minimal eigenvalue of ${\rm  A}$ is positive and is
separated from zero uniformly over $\e,x,t,\o$. By the definitions,
it follows that $\|{\rm  A}_{\e}-{\rm  A}\|_{\W_{\infty}^0}\to 0$.
Since the minimal eigenvalue of a matrix depends continuously of its
coefficients, it follows that the minimal eigenvalue of ${\rm
A}_{\e}$ is positive and is separated from zero uniformly over
$\e,x,t,\o$. Hence Condition \ref{condA}(i) implies (a).
\par
Let us show that Condition \ref{condA}(ii) implies (b).
 Let $
R\defi \|b\|_{\W_{\infty}^1}$, and let
 $\g$ be the supremum  over $x,t,\o$ of
the maximal eigenvalue of $b(x,t,\o)$. It suffices to show that,
without a loss of generality, we can assume that \baa
\sup_{\e}\|b_{\e}\|_{\W_{\infty}^1}\le n\g+2R+1. \label{best}\eaa
\par Suppose that (\ref{best}) does not hold,
i.e., that there exists some $M$ such that for  $\e=M^{-1}$ and some
$t_k=t_k(M)$ there exists $\G\subset\O$ such that $\G\in\F_{t_k}$,
$\P(\G)>0$,  \baaa b_\e(\cdot,t,\o)\equiv b_{\e}(t_k,x,\o),\quad
t\in[t_k,t_{k+1}),\qquad \|b_\e(\cdot,t_k,\o)\|_{W_{\infty}^1(D)}
>n\g+2R+1\quad\hbox{iff}\quad \o\in\G.\eaaa In this case, one
can replace $b_{\e}(\cdot,t)|_{t\in[t_k,t_{k+1})}$, by \baaa \ww
b_{\e}(x,t,\o)=b_\e(x,t,\o)\Ind_{\O\backslash \G}(\o)+ \g
I_n\Ind_{\G}(\o),\quad t\in[t_k,t_{k+1}), \eaaa where $\Ind$ is the
indicator function, and where $I_n$ is the unit matrix in $\R^n$.
Obviously, Condition \ref{condA} is satisfied for $\ww b_\e$
replacing  $ b_\e$, with the same $\d_1>0$ for all $\e$. In
addition, we have that \baaa &&\|\ww
b^{(M)}-b\|_{W^1_{\infty}(D)}\le \|\ww
b^{(M)}|_{W^1_{\infty}(D)}+\|b\|_{W^1_{\infty}(D)}\le n\g+R,\quad
\o\in \G,
 \\
&&\| b^{(M)}-b\|_{W^1_{\infty}(D)}\ge \|\ww
b^{(M)}\|_{W^1_{\infty}(D)}-\|b\|_{W^1_{\infty}(D)}\ge
n\g+2R-R=n\g+R,\quad \o\in\G. \eaaa It follows that Condition
\ref{condTheta} holds for the new selection $\ww b_{\e}$. This
completes the proof of Proposition \ref{proptech}. $\Box$
\vspace{4mm}
\par
Further, it follows from  Proposition \ref{proptech} and from the
definitions that \baa
\sup_{x,t,\o,\e}\Bigl(|b_{\e}(x,t,\o)|+\Bigr|\frac{\p b_{\e}}{\p
x}(x,t,\o)\Bigl|+|f_{\e}(x,t,\o)|+ \Bigr|\frac{\p f_{\e}}{\p
x}(x,t,\o)\Bigl| +|\lambda_{\e}(x,t,\o)|\Bigr)<+\infty.
\label{limits1} \eaa Let us consider a subsequence $\e=\e_i\to 0$
such that \baa b_{\e} \to b,\quad  f_{\e} \to f,\quad
\lambda_{\e}\to \lambda,\quad \frac{\p b_{\e}}{\p x} \to \frac{\p
b_{\e}}{\p x},\quad \frac{\p f_{\e}}{\p x} \to \frac{\p f}{\p
x}\quad \mbox{in}\ X^0\quad \hbox{and a.e.}\quad \hbox{as}\quad
\e\to 0.\quad\label{limits} \eaa
  \par
 Let $p_\e\defi L_\e^*\xi+(\I_T\L_\e)^*\Psi$,
 $\chi_{i\e}\defi \M_\e^*\xi+(\I_T\M_{i\e})^*\Psi$, and let
 $p\defi L^*\xi+(\I_T\L)^*\Psi$,
 $\chi_{i}\defi \M_i^*\xi+(\I_T\M_{i})^*\Psi$.
 The operators $\LU_{\e}^*: X^{-1}\to Y^1$, etc, are defined
similarly to ${{\LU}^*: X^{-1}\to Y^1}$, etc., with substituting
$(b,f,\lambda)=(b_{\e},f_{\e},\lambda_{\e})$.
\par
By Lemma \ref{lemma<1},  the sequences $\{p_{\e}\}$ and
$\{\chi_{i\e}\}$ belong to the closed balls in the spaces $X^2$ and
$X^1$ respectively  with the centers at the zero and with  the
radius $c(\|\xi\|_{X^0}+\|\Psi\|_{Z_T^1})$, where $c=c({\cal P})>0$
does not depend on $\e$. The balls mentioned are closed, concave,
and bounded. It follows that these balls are weakly closed and
weakly compact in the reflexible Banach spaces $X^2$ and $X^1$
respectively. It follows that the sequences $\{p_{\e}\}$ and
$\{\chi_{i\e}\}$ has subsequences with weak limits $\ww p$ and $\ww
\chi_i$, in the corresponding balls, i.e., $$ \|\ww
p\|_{X^2}+\sum_{i=1}^N\|\ww\chi_i\|_{X^1}\le
c(\|\xi\|_{X^0}+\|\Psi\|_{Z_T^1}).$$
\par
 Assume that we can show that $\ww p_\e\to p$  weakly in $X^2$ and $\chi_{i\e}\to \chi_i$
 weakly in $X^1$ for all $i$.
It follows that $\ww p=p$ and $\ww \chi_i=\chi_i$ and \be
\|p\|_{X^2}+ \sum_{i=1}^N\|\chi_i\|_{X^1}\le
c(\|\xi\|_{X^0}+\|\Psi\|_{Z_T^1}). \label{ineq} \ee
 It follows that
$$\Bigl\|\sum_{i=1}^NB_i^*\chi_i\Bigr\|_{X^0}\le c_1
(\|\xi\|_{X^0}+\|\Psi\|_{Z_T^1}),$$ where $c_1=c_1({\cal P})$ is a
constant. Hence  $g\defi \xi+\sum_{i=1}^NB_i^*\chi_i$ is such that
\baaa\|g\|_{X^0}\le
 c_2
(\|\xi\|_{X^0}+\|\Psi\|_{Z_T^1}), \label{gest}\eaaa where
$c_2=c_2({\cal P})$ is a constant.
 Remind that, by Lemma
\ref{prop4.3} and Corollary \ref{corrgxi}, \baaa g=(I-P^*)^{-1}\xi+
(I-P^*)^{-1}P_0^*\Psi, \quad p=\Q_0^*g+(\I_T\Q_0)^*\Psi. \eaaa  By
Lemma \ref{lemmaQ}, it follows that $p\in Y^2$
 and $$ \|p\|_{Y^2}+\Bigl\|\sum_{i=1}^NB_i^*\chi_i\Bigr\|_{X^0}\le
 c_3
(\|g\|_{X^0}+\|\Psi\|_{Z_T^1}),$$ where $c_3=c_3({\cal P})$. By
(\ref{gest}), it follows that $p\in Y^2$
 and $$ \|p\|_{Y^2}+\Bigl\|\sum_{i=1}^NB_i^*\chi_i\Bigr\|_{X^0}\le
 c_4
(\|\xi\|_{X^0}+\|\Psi\|_{Z_T^1}),$$ where $c_4=c_4({\cal P})$ is a
constant. Then the proof of Theorem \ref{Th2FE} follows provided
that the weak convergence of the sequence $\{\chi_{i\e}\}$ to
$\chi_i$ is established.
\par
 Therefore, it suffices prove this weak convergence, i.e., it suffices to  show that
\baa &&I_{\e}\defi (p_{\e}-p,h)_{X^0} \to 0 \quad \mbox{as} \quad\e
\to 0\qquad \forall h\in X^0,  \label{5.7a}\\
&&J_{\e}\defi (\chi_{i\e}-\chi_i,h)_{X^0} \to 0 \qquad \mbox{as}
\quad\e \to 0\quad \forall h\in X^0, \ i\in\{1,...,N\}. \label{5.7b}
 \eaa
\par
Let us show that  (\ref{5.7a}) holds. Set $u_{\e}\defi \LU_{i\e}\,h$
and $u\defi \LU_i\,h$, where the operators $\LU_{i\e}: X^{0}\to Y^1$
are defined similarly to the operators $\LU_i: X^{0}\to Y^1$ with
substituting $(b,f,\lambda)=(b_{\e},f_{\e},\lambda_{\e})$. By the
definitions of the corresponding adjoint operators, \baaa
I_{\e}&=&(\LU_{i\e}^*\xi-\LU_i^* \xi
,h)_{X^0}+((\I_T\LU_{i\e})^*\Psi-(\I_T\LU_i)^* \Psi ,h)_{X^0}
\\&=&(\xi,u_{\e}-u)_{X^0}+\Big(\Psi,u_{\e}(\cdot,T)-u(\cdot,T)\Big)_{Z_T^0}.
\eaaa
\par
Let the operators $\A_{\e}$ be defined similarly to $\A$ with
substituting $(b,f,\lambda)=(b_{\e},f_{\e},\lambda_{\e})$. By the
definitions, it follows that there exist functions $\w f_\e(x,t,\o):
\R^n\times\R_+\times\O\to \R^n$, and $\w \lambda_\e (x,t,\o):
\R^n\times\R_+\times\O\to \R$, such that \baaa
 \sup_{\e>0}\esssup_{x,t,\o}\Bigl(|\w f_\e(x,t,\o)|+|\w \lambda_\e (x,t,\o)|\Bigr)<+\infty,
 \eaaa
and that $\A_{\e} u-\A u$ is represented as
 \baaa
\A_{\e} u-\A u =  \sum_{i,j=1}^n \frac{\p }{\p x_i }\Bigl(
[b_{ij\e}-b_{ij}]\,\frac{\p u}{\p x_j}\,\Bigr) + \sum_{i=1}^n \w
f_{i\e}\frac{\p u}{\p x_i }\,+\w \lambda_\e\,u. \label{AA}\eaaa By
(\ref{limits1})-(\ref{limits}), it follows that \baa \w f_{\e} \to
0\quad \mbox{and}\quad \w\lambda_{\e}\to 0\quad \mbox{in}\ X^0\quad
\hbox{and a.e.}.\label{limits3} \eaa
 The function $U_{\e}\defi u_{\e}-u$ is the solution in $Q$ of
the boundary value problem \baaa &&d_t U_{\e}
=(\A_{\e}U_{\e}+F_{\e}(u))\,dt +\sum_{i=1}^NB_iU_{\e}\,dw_i(t),
\\
&& U_{\e}(x,0)=0,\quad U_{\e}(x,t)|_{x\in\p D}=0, \eaaa and where
the linear operator $F_{\e}(\cdot)$ is defined as
\begin{eqnarray*}
F_{\e}(u)\defi r_{\e}(u)+q_{\e}(u), \qquad
r_{\e}(u)\defi\sum_{i,j=1}^n\frac{\p }{\p x_i}\Bigl([\w b_{\e
ij}-b_{ij}]\frac{\p u}{\p x_j}\Bigr), \qquad q_{\e}(u)\defi\frac{\p
u}{\p x} \,\w f_{\e}+ \w \lambda_{\e}u.
\end{eqnarray*}
 Here
$\w b_{\e ij}$ are the components of the matrix $\w b_{\e}$.
 By Lemma \ref{lemma1}, it follows that $$\|
U_{\e}\|_{Y^1}\le C \|F_{\e}(u)\|_{X^{-1}}, $$ for a constant
$C_1=C_1({\cal P})$. It follows that there exists a constant
$C=C({\cal P})>0$ such that \baaa | I_{\e}| \le C\|
U_{\e}\|_{Y^1}(\|\xi\|_{X^{-1}}+\|\Psi\|_{Z_0^T}) \le C
\|F_{\e}(u)\|_{X^{-1}}(\|\xi\|_{X^{0}}+\|\Psi\|_{Z_1^T}).\label{Iest}\eaaa
\par We have that \baaa
\|r_{\e}(u)\|_{X^{-1}}^2=\E \int_0^T\Bigl\|\sum_{i,j=1}^n\frac{\p
}{\p x_i}\Bigl([\w b_{\e ij}-b_{ij}]\frac{\p u}{\p
x_j}\Bigr)\Bigr\|_{H^{-1}}^2dt \le C_1\sum_{i,j=1}^n\E
\int_0^T\Bigl\|[\w b_{\e ij}-b_{ij}]\frac{\p u}{\p
x_j}\Bigr\|_{H^{0}}^2dt,  \eaaa for a constant $C=C(n)$. The
functions $b_{\e}$ and $b$ are bounded, hence $$ \Bigl| [\w b_{\e
ij}-b_{ij}]\frac{\p u}{\p x_j} \Bigr| \le C_1\Bigl| \frac{\p u}{\p
x}(x,t,\o)\Bigr|
$$ for a constant $C_1=C_1({\cal P})$. We have that $u \in X^1$.  By the Lebesgue's Dominated
Convergence Theorem,  it follows that $\Bigl\|[\w b_{\e
ij}-b_{ij}]\frac{\p u}{\p x_j}\Bigr\|_{X^{0}}\to 0$. Hence
$\|r_{\e}(u)\|_{X^{-1}}\to 0$.
\par
Further, the functions $\w f_{\e}$ and $\w \lambda_{\e}$ are
bounded, hence
$$ | q_{\e}(u)(x,t,\o)| \le C_1\bigg(\bigg| \frac{\p u}{\p
x}(x,t,\o)\bigg|+ |u(x,t,\o)|\bigg)
$$ for a constant $C_2=C_2({\cal P})
>0$.
 By the Lebesgue's Dominated
Convergence Theorem again,  it follows that ${\|
q_{\e}(u)\|_{X^0}\to 0}$. Therefore, we obtain that ${\|
U_{\e}(u)\|_{X^0}\to 0}$. By (\ref{Iest}),  it follows that
(\ref{5.7a}) holds.
\par
Let us show that  (\ref{5.7b}) holds. Set $v_{\e}\defi \M_{i\e}\,h$
and $v\defi \M_i\,h$, where the operators $\M_{i\e}: X^{0}\to Y^1$
are defined similarly to the operators $\M_i: X^{0}\to Y^1$ with
substituting $(b,f,\lambda)=(b_{\e},f_{\e},\lambda_{\e})$. By the
definitions of the corresponding adjoint operators,  \baaa
J_{\e}&=&(\M_{i\e}^*\xi-\M_i^* \xi
,h)_{X^0}+((\I_T\M_{i\e})^*\Psi-(\I_T\M_i)^* \Psi ,h)_{X^0}
\\&=&(\xi,v_{\e}-v)_{X^0}+\Big(\Psi,v_{\e}(\cdot,T)-v(\cdot,T)\Big)_{Z_T^0}.
\eaaa  The function $V_{\e}\defi v_{\e}-v$ is the solution in $Q$ of
the boundary value problem \baaa &&d_t V_{\e}
=(\A_{\e}V_{\e}+F_{\e}(v))\,dt +\sum_{i=1}^NB_iV_{\e}\,dw_i(t),
\\
&& V_{\e}(x,0)=0,\quad V_{\e}(x,t)|_{x\in\p D}=0, \eaaa where the
operator $F_{\e}(\cdot)$ is defined above. The remaining part of the
proof of (\ref{5.7b}) repeats the proof of (\ref{5.7a}).
 This completes the proof of Theorem \ref{Th2FE}.
$\Box$
\section{The proof of Theorems \ref{propNewOld}-\ref{propON2} and \ref{ThRobust}}}
{\it Proof of Theorem \ref{propNewOld}.} Assume that Condition
\ref{condA} holds.  Let \baaa
S_N\defi\Bigl\{\a=(\a_1,...,\a_N)^\top\in\R^N:\quad
|\a|=\Bigl(\sum_{i=1}^N|\a_i|^2\Bigr)^{1/2}\le 1\Bigr\}.\eaaa Let
$y\in\R^n$ be fixed and let $y_i=y_i(\a)\defi \a_i y$, $\a\in S_N$.
Let $y_i\defi \a_i y$ and $z_i=z_i(y)=\b_i^\top y$,
$z=z(y)=(z_1,...,z_N)^\top$. By Condition \ref{condA}, \baaa
 \label{Main1--} \sum_{i=1}^Ny_i^\top  b
\,y_i\ge\frac{1}{2}\left(\sum_{i=1}^Ny_i^\top\b_i\right)^2 +
\d_1\sum_{i=1}^N|y_i|^2 \eaaa for all $\a\in S_N$, $(x,t)\in D\times
[0,T]$ and $\o\in\O$. Hence \baaa y^\top b
\,y=\sum_{i=1}^N\a_i^2y^\top b \,y\ge
\frac{1}{2}\left(\sum_{i=1}^N\a_iy^\top\b_i\right)^2 +
\d_1\sum_{i=1}^N\a_i^2|y|^2 =
\frac{1}{2}\left(\sum_{i=1}^N\a_iz_i(y)\right)^2 +
\d_1|y|^2\sum_{i=1}^N\a_i^2\\=\frac{1}{2}(\a^\top z(y))^2 +
\d_1|y|^2 \eaaa for any $\a\in S_N$. Hence \baaa
 y^\top  b
\,y\ge \sup_{\a\in S_N}\frac{1}{2}\left(\a^\top z(y)\right)^2 +
\d_1|y|^2= \frac{1}{2}|z(y)|^2 + \d_1|y|^2.  \eaaa On the other
hand, \baaa |z(y)|^2 = \sum_{i=1}^N |z_i(y)|^2=\sum_{i=1}^N
|y^\top\b_i|^2.\eaaa Hence \baaa
  y^\top  b
\,y\ge \frac{1}{2}\sum_{i=1}^N |y^\top\b_i|^2 + \d_1|y|^2. \eaaa
Hence Condition \ref{cond3.1.A} holds with $\d=\d_1$.
  $\Box$
\par
{\it Proof of Theorem \ref{propON1}.} We have that $2b=\g+R$, where
$\g=\sum_{i=1}^n\b_i^2$ and $R=R(x,t,\o)\ge 2\d$.  Let $D\defi
BB^\top =\{\b_i\b_j\}_{i,j=1}^N$, where $B\defi
(\b_1,...,\b_N)^\top$. It suffices to show that there exists
$\d_1>0$ such that \baa \g(x,t,\o) I_N-D(x,t,\o)\ge 0\eaa for all
$x,t,\o$,  where $I_N$ is the unit matrix in $\R^{N\times N}$. Let
$\lambda=\lambda(x,t,\o)$ be the minimal eigenvalue of the matrix
$\g(x,t,\o) I_N-D(x,t,\o)$. It suffices to show that $\lambda \ge
0$. Let $z=z(x,t,\o)$ be a corresponding eigenvector such that
$|z|=|B|\neq 0$ (for the trivial case $|B|=0$, we have immediately
that $\lambda =0$). We have that  $z=cB+B'$, where $c\in [-1,1]$ and
$B'=B'(x,t,\o)$ is a vector such that $B^\top B'=0$. By the
definitions, we have that $\g=|B|^2$ and  \baaa \lambda z =(\g
I_N-D)z=(\g I_N-BB^\top) (c B+B')=\g (cB +B')-c|B|^2B\\= \g cB+\g
B'-c\g B=\g B'. \eaaa Hence $\lambda (cB+B')=\g B'.$  It follows
that either $B'\neq 0$, $c=0$,  and $\lambda = \g\ge 0$, or
$\lambda=0$ and $B'=0$. This completes the proof. $\Box$.
\par
{\it Proof of Theorem \ref{propON2}.} By H\"older inequality, we
have that \baaa \left(\sum_{i=1}^{N_0}y_i^\top\b_i \right)^2\le
N_0\sum_{i=1}^{N_0}\left(y_i^\top\b_i \right)^2. \eaaa Hence
 \baaa \sum_{i=1}^Ny_i^\top  b
\,y_i-\frac{1}{2}\left(\sum_{i=1}^Ny_i^\top\b_i\right)^2 =
\sum_{i=1}^{N_0}y_i^\top  b
\,y_i-\frac{1}{2}\left(\sum_{i=1}^{N_0}y_i^\top\b_i\right)^2 \ge
\sum_{i=1}^{N_0}y_i^\top  b
\,y_i-\frac{N_0}{2}\sum_{i=1}^{N_0}\left(y_i^\top\b_i\right)^2\nonumber\\\ge
\d_2\sum_{i=1}^N|y_i|^2.\hphantom{xxxx} \eaaa This completes the
proof. $\Box$.
\par
{\it Proof of Theorem \ref{ThRobust}}. Let  \baaa p\defi
p^{(1)}-p^{(2)},\quad \chi_i\defi \chi_i^{(1)}-\chi_i^{(2)}, \eaaa
and let $\A^{(k)*}$, $B^{(k)*}_i$ be the corresponding operators
(\ref{AB*}), $k=1,2$. We have that \baaa \label{parabR2}
&&d_tp+(\A^{(1)*} p+ \psi)\,dt +\sum_{i=1}^N B_i^{(1)}\chi_i dt
+\psi=\chi_i\,dw_i(t), \quad t\le T,
\\
&&p(x,0,\o)=\Psi^{(1)}(x,\o)-\Psi^{(2)}(x,\o), \qquad
p(x,t,\o)\,|_{x\in \p D}=0.
\label{parabR22}
\eaaa Here
$$
\psi\defi \xi^{(1)}-\xi^{(2)}+ \A^{(1)*}p^{(2)}-
\A^{(2)*}p^{(2)}+\sum_{i=1}^N (B_i^{(1)*}\chi_i^{(2)}-
B_i^{(2)*}\chi_i^{(2)}).
$$
By Theorem \ref{Th2FE}, it follows that  there exists a constant
$C_0=C_0({\cal P}^{(1)})$  such that \be\label{UU}\| p \|_{{Y}^2}+
\sum_{i=1}^N\|\chi_i\|_{X^1}
 \le   C_0
 \biggl(  \|\psi\|_{X^0}
+\|\Psi^{(1)}-\Psi^{(2)}  \|_{Z_0^1}\biggr). \ee Further,  we have $
\psi=\sum_{m=0,1,2}\psi_m, $ where \baaa
&&\psi_0=\xi^{(1)}-\xi^{(2)},\\ &&\psi_1\defi
\sum_{i,j=1}^n[b_{ij}^{(1)}-b_{ij}^{(2)}] \frac{\p^2p^{(2)}}{\p
x_i\p x_j} +\sum_{i=1}^n[f_{i}^{(1)}-f_{i}^{(2)}] \frac{\p
p^{(2)}}{\p x_i}+[\lambda^{(1)}-\lambda^{(2)}]p^{(2)},
 \\ &&\psi_2= \sum_{i=1}^N
 \left(\sum_{i=1}^n[\b_{i}^{(1)}-\b_{i}^{(2)}]
\frac{\p \chi_i^{(2)}}{\p x_i}+[\w\b^{(1)}-\w\b^{(2)}]
\chi_i^{(2)}\right). \eaaa
 Clearly,
$$
\|\psi_0\|_{X^0} \le M, \qquad \|\psi_1\|_{X^0} \le C
M\|p^{(2)}\|_{Y_2},\qquad |\psi_2\|_{X^0} \le C M
+\sum_{i=1}^N\|\chi^{(2)}\|_{X^1},
$$
where $C=C(n)$ is a constant.  Finally, we obtain
$$
\|\psi\|_{X^0} \le
C_1M(\|p^{(2)}\|_{Y_2}+\sum_{i=1}^N\|\chi^{(2)}_i\|_{X^1} +1),
$$
where $C_1=C_1(n)$ is a constant. By (\ref{UU}), the desired
estimate follows. This completes the proof. $\Box$
\subsection*{Acknowledgment}  This work  was supported
by NSERC
grant of Canada 341796-2008 to the author.

\section*{References}
$\hphantom{XX}$Al\'os, E., Le\'on, J.A., Nualart, D. (1999).
 Stochastic heat equation with random coefficients
 {\it
Probability Theory and Related Fields} {\bf 115}, 1, 41-94.
\par
Bally, V., Gyongy, I., Pardoux, E. (1994). White noise driven
parabolic SPDEs with measurable drift. {\it Journal of Functional
Analysis} {\bf 120}, 484 - 510.
\par
Chojnowska-Michalik, A., and Goldys, B. (1965). {Existence,
uniqueness and invariant measures for stochastic semilinear
equations in Hilbert spaces},  {\it Probability Theory and Related
Fields},  {\bf 102}, No. 3, 331--356.
\par
Confortola, F. (2007).  Dissipative backward stochastic
differential equations with locally Lipschitz nonlinearity, {\it
Stochastic Processes and their Applications} {\bf 117}, Issue 5,
613-628.
\par
Da Prato, G., and Tubaro, L. (1996). { Fully nonlinear stochastic
partial differential equations}, {\it SIAM Journal on Mathematical
Analysis} {\bf 27}, No. 1, 40--55.\
\par
Dokuchaev, N.G. (1992). { Boundary value problems for functionals
of
 Ito processes,} {\it Theory of Probability and its Applications}
 {\bf 36 }, 459-476.
\par
 Dokuchaev, N.G. (1995). Probability distributions  of  Ito's
processes: estimations for density functions and for conditional
expectations of integral functionals. {\it Theory of Probability and
its Applications} {\bf 39} (4),  662-670.
\par Dokuchaev, N.G. (2003). Nonlinear
parabolic Ito's equations and duality approach, {\it Theory of
Probability and its Applications} {\bf 48} (1), 45-62.
\par
Dokuchaev, N.G. (2005). Parabolic Ito equations and second
fundamental inequality.  {\it Stochastics} {\bf 77}, iss. 4.,
349-370.
\par
Dokuchaev, N.G. (2008). Universal estimate of the  gradient for
parabolic equations. {\it Journal of Physics A: Mathematical and
Theoretical} {\bf 41}, No. 21,  215202 (12pp).
\par
Dokuchaev, N. (2010a). Duality and semi-group property for backward
parabolic Ito equations. {\em Random Operators and Stochastic
Equations. } {\bf 18}, 51-72.     (See also working paper at
arXiv:0708.2497v2 [math.PR] at http://lanl.arxiv.org/abs/0708.2497).
\par
Dokuchaev, N. (2010b). Representation of functionals  of Ito
processes in bounded domains. {\em Accepted to Stochastics}. (See
also working paper at    arXiv:math/0606601v2 [math.PR]
http://lanl.arxiv.org/abs/math/0606601).
\par
Gy\"ongy, I. (1998). Existence and
uniqueness results for semilinear stochastic partial differential
equations. {\it Stochastic Processes and their Applications} {\bf
73} (2), 271-299.
\par
Hu, Y., Peng, S. (1991). Adapted solution of a backward semilinear
stochastic evolution equation. {\it Stochastic Anal. Appl.} {\bf
9}, 445-459.
\par
 Hu, Y., Ma, J., and Yong, J. (2002).  On Semi-linear Degenerate
Backward Stochastic Partial Differential Equations. {\it
Probability Theory and Related Fields}, {\bf 123}, No. 3, 381-411.
\par
Kim, Kyeong-Hun (2004). On stochastic partial di!erential
equations with variable coefficients in $C^1$-domains. {\it
Stochastic Processes and their Applications} {\bf 112}, 261--283.
\par Krylov, N. V. (1999). An
analytic approach to SPDEs. {\it Stochastic partial differential
equations: six perspectives}, 185--242, Math. Surveys Monogr., 64,
Amer. Math. Soc., Providence, RI.
\par Ladyzhenskaia, O.A. (1985).
{\it The Boundary Value Problems of Mathematical Physics}. New York:
Springer-Verlag.
\par
Ladyzenskaya, O.A., Solonnikov, V.A., and Ural'ceva, N.N. (1968).
{\it Linear and quasi--linear equations of parabolic type.}
Providence, R.I.: American Mathematical Society.
\par Ma, J., and Yong, J. (1999),  On Linear Degenerate
Backward Stochastic PDE's. {\it Probability Theory and Related
Fields}, {\bf 113}, 135-170.
\par
Maslowski, B. (1995). { Stability of semilinear equations with
boundary and pointwise noise}, {\it Ann. Scuola Norm. Sup. Pisa
Cl. Sci.} (4), {\bf  22}, No. 1, 55--93.
\par
Pardoux, E., S. Peng, S. (1990). Adapted solution of a backward
stochastic differential equation. {\it System \& Control Letters}
{\bf 14}, 55-61.
\par
Pardoux, E. (1993).
 Stochastic partial differential equations, a review, {\it Bull. Sc. Math.}
 {\bf 117}, 29-47.
 \par
Pardoux, E., A. Rascanu,A. (1998). Backward stochastic
differential equations with subdifferential operators and related
variational inequalities, {\it Stochastic Process. Appl.} {\bf 76}
(2), 191-215.
\par
Rozovskii, B.L. (1990). {\it Stochastic Evolution Systems; Linear
Theory and Applications to Non-Linear Filtering.} Kluwer Academic
Publishers. Dordrecht-Boston-London.
\par
Tessitore, G. (1996). Existence, uniqueness and space regularity
of the adapted solutions of a backward SPDE. Stochastic Anal.
Appl. 14 461-486.
\par
Walsh, J.B. (1986). An introduction to stochastic partial
differential equations, {\it Ecole d'Et\'e de Prob. de St.} Flour
XIV, 1984, Lect. Notes in Math 1180, Springer Verlag.
\par
Yong, J., and Zhou, X.Y. (1999). { Stochastic controls:
Hamiltonian systems and HJB equations}. New York: Springer-Verlag.
\par
Yosida, K. (1965).   {\it Functional Analysis}. Springer-Verlag.
Berlin, Gottingen, Heidelberg.
\par
 Zhou, X.Y. (1992). { A duality analysis on stochastic partial
differential equations}, {\it Journal of Functional Analysis} {\bf
103}, No. 2, 275--293.
\end{document}